\documentclass[11pt]{amsart}

\usepackage{amsmath}
\usepackage{amscd}
\usepackage{amssymb}
\usepackage{latexsym}
\usepackage{stmaryrd} 
\usepackage{color}
\usepackage{mathrsfs}

\usepackage[arrow,curve,matrix,tips,2cell]{xy}
  \SelectTips{eu}{10} \UseTips
  \UseAllTwocells
\newtheorem{theorem}{Theorem}[section]
\newtheorem*{theorem*}{Theorem}
\newtheorem{lemma}[theorem]{Lemma}
\newtheorem*{lemma*}{Lemma}
\newtheorem{corollary}[theorem]{Corollary}
\newtheorem*{corollary*}{Corollary}
\newtheorem{proposition}[theorem]{Proposition}

\newtheorem{remark}[theorem]{Remark}
\newtheorem{definition}[theorem]{Definition}

\newtheorem{question}[theorem]{Question}
\newtheorem{example}[theorem]{Example}
\newcommand{\bgl}{\begin{equation}} 
\newcommand{\egl}{\end{equation}}
\newcommand{\bgloz}{\begin{equation*}} 
\newcommand{\egloz}{\end{equation*}}
\newcommand{\bgln}{\begin{eqnarray}} 
\newcommand{\egln}{\end{eqnarray}}
\newcommand{\bglnoz}{\begin{eqnarray*}} 
\newcommand{\eglnoz}{\end{eqnarray*}}
\newcommand{\btheo}{\begin{theorem}}
\newcommand{\etheo}{\end{theorem}}
\newcommand{\btheooz}{\begin{theorem*}}
\newcommand{\etheooz}{\end{theorem*}}
\newcommand{\blemma}{\begin{lemma}}
\newcommand{\elemma}{\end{lemma}}
\newcommand{\blemmaoz}{\begin{lemma*}}
\newcommand{\elemmaoz}{\end{lemma*}}
\newcommand{\bproof}{\begin{proof}}
\newcommand{\eproof}{\end{proof}}
\newcommand{\bbew}{\begin{beweis}}
\newcommand{\ebew}{\end{beweis}}
\newcommand{\bremark}{\begin{remark}\em}
\newcommand{\eremark}{\end{remark}}
\newcommand{\bdefin}{\begin{definition}}
\newcommand{\edefin}{\end{definition}}
\newcommand{\bprop}{\begin{proposition}}
\newcommand{\eprop}{\end{proposition}}
\newcommand{\bcor}{\begin{corollary}}
\newcommand{\ecor}{\end{corollary}}
\newcommand{\bcoroz}{\begin{corollary*}}
\newcommand{\ecoroz}{\end{corollary*}}
\newcommand{\bfa}{\begin{cases}} 
\newcommand{\efa}{\end{cases}}
\newcommand{\bquestion}{\begin{question}\em}
\newcommand{\equestion}{\end{question}}
\newcommand{\bexample}{\begin{example}\em}
\newcommand{\eexample}{\end{example}}
\newcommand{\cF}{\mathcal F}

\newcommand{\cQ}{\mathcal Q}
\newcommand{\cR}{\mathcal R}
\newcommand{\cS}{\mathcal S}
\newcommand{\cT}{\mathcal T}

\newcommand{\cV}{\mathcal V}

\def\Nz{\mathbb{N}}

\def\Zz{\mathbb{Z}}

\def\1z{\mathbb{1}}
\newcommand{\sR}{\mathscr{R}}
\newcommand{\sS}{\mathscr{S}}
\newcommand{\ti}{\tilde}
\newcommand{\wt}{\widetilde}

\newcommand{\lori}{\longrightarrow}

\newcommand{\ma}{\mapsto} 
\newcommand\onto{\twoheadrightarrow} 
\newcommand\into{\hookrightarrow} 

\newcommand{\ve}{\varepsilon}

\def\SEMI{\mbox{$\times\kern-2pt\vrule height5pt width.6pt \kern3pt $}}

\newcommand{\id}{{\rm id}}

\renewcommand{\ker}{{\rm ker}\,}

\newcommand{\reg}{^\times} 
\newcommand{\lspan}{{\rm span}} 
\newcommand{\abs}[1]{\lvert#1\rvert} 
\newcommand{\defeq}{\mathrel{:=}} 
\newcommand{\eqdef}{\mathrel{=:}} 

\newcommand{\dop}{\text{: }} 
\newcommand{\fa}{\text{ for all }} 
\newcommand{\ilim}{\varinjlim} 

\newcommand{\Idem}{{\rm Idem}\,}

\newcommand{\lge}{\left\{} 
\newcommand{\rge}{\right\}} 
\newcommand{\lsp}{\left\langle} 
\newcommand{\rsp}{\right\rangle} 
\newcommand{\gekl}[1]{\lge #1 \rge} 
\newcommand{\spkl}[1]{\lsp #1 \rsp} 
\newcommand{\menge}[2]{\gekl{ #1 \dop #2 }} 
\begin{document}

\title[Independent resolutions I: Algebraic case]{Independent resolutions for totally disconnected dynamical systems I: Algebraic case}

\author{Xin Li}
\author{Magnus Dahler Norling}

\address{Xin Li, School of Mathematical Sciences, Queen Mary University of London, Mile End Road, London E1 4NS, UK}
\email{xin.li@qmul.ac.uk}
\address{Magnus Dahler Norling, Institute of Mathematics, University of Oslo, P.b. 1053 Blindern, 0316 Oslo, Norway}
\email{magnudn@math.uio.no}


\thanks{\scriptsize{The research of the first named author was partially supported by the ERC through AdG 267079.}}

\begin{abstract}
This is the first out of two papers on independent resolutions for totally disconnected dynamical systems. In the present paper, we discuss independent resolutions from an algebraic point of view. We also present applications to group homology and cohomology. This first paper sets the stage for our second paper, where we explain how to use independent resolutions in K-theory computations for crossed products attached to totally disconnected dynamical systems.
\end{abstract}

\maketitle


\setlength{\parindent}{0pt} \setlength{\parskip}{0.5cm}

\section{Introduction}

This paper is about algebraic independent resolutions. Our original motivation to introduce independent resolutions came from K-theory for C*-algebras. Recently, in \cite{{C-E-L1},{C-E-L2}}, a new method was developed to compute K-theory for crossed products attached to certain dynamical systems on totally disconnected spaces. This method only works for dynamical systems admitting a so-called invariant regular basis. By \cite[Definition~2.9]{C-E-L2} a regular basis for the totally disconnected Hausdorff space $\Omega$ is a collection $\cV$ of non-empty compact open subsets of $\Omega$ such that $\cV$ is closed under finite intersections as long as they are non-empty, any compact open subset of $\Omega$ can be constructed from elements of $\cV$ by a finite number of intersections, unions and set differences (i.e. $\cV$ generates the collection of all compact open subsets of $\Omega$ as a Boolean algebra), and proper finite unions of sets in $\cV$ are not in $\cV$ again. If a group acts continuously on $\Omega$, then a regular basis $\cV$ for $\Omega$ is called invariant if the group action maps $\cV$ into itself. The existence of an invariant regular basis occurs in natural examples, for instance if our dynamical systems come from certain subsemigroups of groups. However, this condition can also be quite restrictive, as observed in \cite[Proposition~3.18]{C-E-L2}. Therefore, a natural task is to extend the method from \cite{C-E-L2} to more general dynamical systems.

The main idea is as follows: Given an arbitrary dynamical system which might not admit an invariant regular basis, we want to find a sequence of dynamical systems, each of which admitting invariant regular bases, such that this sequence leads to a resolution of our original dynamical system. In terms of C*-algebras, this means that the crossed products attached to our dynamical systems fit into a long exact sequence whose last non-zero entry is the crossed product of the dynamical system we started with. Since all the dynamical systems in our sequence admit invariant regular bases, we can apply the K-theoretic method from \cite{C-E-L2} to the corresponding crossed products. Together with our long exact sequence, this leads to K-theoretic computations for the crossed product of our original dynamical system.

In the process of developing this idea, we realized that our methods actually work in a purely algebraic context. Namely, since we are looking at totally disconnected spaces, it makes sense to replace complex-valued continuous functions by integer-valued ones. The algebra of such functions is just a $\Zz$-algebra. Given an arbitrary dynamical system on a totally disconnected space, we are able to produce a sequence of totally disconnected dynamical systems which all admit invariant regular bases and whose $\Zz$-algebras (in the sense above) give rise to a resolution of the $\Zz$-algebra of our original dynamical system (or rather of its underlying space). The point is that everything works equivariantly with respect to the group actions. We call such a resolution an algebraic independent resolution. This notion leads to C*-algebraic independent resolutions for the purpose of computing K-theory, in the sense as explained above. But not only this; already on a purely algebraic level, our notion of algebraic independent resolutions has interesting applications. To give a concrete example, we explain how to compute group homology and cohomology using algebraic independent resolutions. Moreover, since we work in the category of $\Zz$-algebras, these resolutions might be interesting for studying homology or cohomology of algebras. Our notion of algebraic independent resolutions might also be interesting from the point of view of dynamical systems. Totally disconnected dynamical systems arise naturally in various contexts, for instance in symbolic dynamics. And given such a system, the minimal length of an algebraic independent resolution is an invariant of the dynamical system which might be worth investigating.

For applications to K-theory for C*-algebras, we refer the reader to our second paper \cite{L-N}, where we discuss graph C*-algebras, C*-algebras of one dimensional tilings, various ideals and quotients of semigroup C*-algebras, and semigroup C*-algebras of semigroups not satisfying the independence condition which could not be treated with the methods in \cite{C-E-L2}.

We are particularly interested in independent resolutions of finite length. These are much easier to treat than infinite ones, especially when it comes to K-theory. Therefore, one of our main goals is to give a criterion for the existence of finite length independent resolutions. The idea is the following: The algebra of integer-valued continuous functions on a totally disconnected space is generated by idempotents. If our dynamical system does not admit an invariant regular basis, then these idempotents - if they are invariant under the group action - must satisfy non-trivial (integral, linear) relations. We give a sufficient condition for the existence of finite length independent resolutions in terms of these relations: If these relations satisfy a certain finiteness condition, then we are able to produce a finite length independent resolution.

\section{Algebraic independent resolutions}
\label{sec-ind-res}

The situation we would like to consider can be viewed from two different angles. From the point of view of dynamical systems, we want to study a dynamical system $\Gamma \curvearrowright \Omega$ with a discrete group $\Gamma$ acting by homeomorphisms on a totally disconnected (locally compact Hausdorff) space $\Omega$. From the point of view of algebra, we are interested in an action of a group $\Gamma$ on a commutative $\Zz$-algebra $D$ generated by idempotents. Since every commutative $\Zz$-algebra generated by idempotents is precisely of the form $C_0(\Omega, \Zz)$ for a totally disconnected space $\Omega$, these two viewpoints are equivalent. In the following, we introduce the notion of an algebraic independent resolution of $\Gamma \curvearrowright C_0(\Omega, \Zz)$, or equivalently, $\Gamma \curvearrowright D$.

Let us start with some notation.
\bdefin
A semilattice $E$ is a commutative idempotent semigroup, i.e., a set with an associative and commutative binary operation in which every element $e$ satisfies $e e = e$.
\edefin
All our semilattices will have a zero element, denoted $0$, which satisfies $0e=0$ for each $e\in E$. On the semilattice $E$ there is a partial order, called the natural partial order, given as follows:
\bdefin
For $e, f \in E$, we say that $e$ is smaller than $f$, or that $f$ dominates $e$, if $e f = e$. We write $e \leq f$ in that case. Also, we write $e \lneq f$ if $e \leq f$ and $e \neq f$.
\edefin

Given a semilattice $E$, we set $E\reg \defeq E \setminus \gekl{0}$. Moreover, the integral semigroup ring $\Zz[E\reg]$ of $E\reg \subseteq E$ is given by the free $\Zz$-module $\oplus_{e \in E\reg} \Zz$ on $E\reg$. We write $\sum_e m_e e$ for the element in $\Zz[E\reg]$ whose entry at $e \in E\reg$ is $m_e \in \Zz$. Multiplication in $\Zz[E\reg]$ is given by
$$(\sum_e m_e e)(\sum_f n_f f) = \sum_{\substack{e, f \\ ef \neq 0}} n_e n_f (ef).$$
Note that if $\Zz[E]$ denotes the usual semigroup ring over $\Zz$ with free $\Zz$-basis $E$, then $\Zz[E\reg]$ is canonically isomorphic to $\Zz[E] / (\Zz \cdot 0)$. $\Zz[E\reg]$ is characterized by the following universal property: Given a $\Zz$-algebra $A$, let $\Idem(A) \defeq \menge{a \in A}{a = a^2}$. Whenever $f: E \to \Idem(A)$ is a semigroup homomorphism satisfying $f(0) = 0$, there exists a unique $\Zz$-algebra homomorphism $\Zz[E\reg] \to A$ sending $e$ to $f(e)$. In particular, a semigroup homomorphism $f: \: E_1 \to E_2$ induces a (uniquely determined) $\Zz$-algebra homomorphism $\Zz[E_1\reg] \to \Zz[E_2\reg]$ sending $e \in E_1\reg$ to $f(e) \in E_2\reg$ if $f(e) \neq 0$ and to $0 \in \Zz[E_2\reg]$ if $f(e) = 0$. In other words, the assignment $E \ma \Zz[E\reg]$ is functorial.

By an action of a group $\Gamma$ on a semilattice $E$ we mean a group homomorphism from $\Gamma$ to the semigroup automorphisms of $E$. A $\Gamma$-semilattice is a semilattice together with a group action of $\Gamma$. Such an action obviously induces an action of $\Gamma$ on $\Zz[E\reg]$.

It turns out that every $\Zz$-algebra of the form $C_0(\Omega, \Zz)$ for a totally disconnected space $\Omega$ is isomorphic to the integral semigroup ring of a suitable semilattice. Namely, by \cite[Proposition~2.12]{C-E-L2}, we can always find a regular basis $\cV$ for $\Omega$ in the sense of \cite[Definition~2.9]{C-E-L2}. Recall that a regular basis is a collection $\cV$ of non-empty compact open subsets of $\Omega$ such that $\cV$ is closed under finite intersections as long as they are non-empty, any compact open subset of $\Omega$ can be constructed from elements of $\cV$ by a finite number of intersections, unions and set differences, and proper finite unions of sets in $\cV$ are not in $\cV$ again. Given such a regular basis $\cV$, $E \defeq \menge{e_V}{V \in \cV} \cup \gekl{0}$ is a semilattice with respect to $e_V e_W \defeq e_{V \cap W}$ if $V \cap W \neq \varnothing$, and $e_V e_W \defeq 0$ if $V \cap W = \varnothing$ (multiplication with $0$ gives $0$, as always). It is clear that if we view $\cV \cup \gekl{\varnothing}$ as a semilattice with respect to intersection, the map $\cV \cup \gekl{\varnothing} \to E, \, V \ma e_V, \, \varnothing \ma 0$ is an isomorphism of semilattices. The same argument as in \cite[Remark~3.22]{C-E-L2}, but for integral semigroup rings in place of C*-algebras of semilattices, yields the isomorphism $\Zz[E\reg] \cong C_0(\Omega, \Zz)$, $e_V \ma 1_V$. Here $1_V$ is the characteristic function of $V$.

Now given a totally disconnected dynamical system $\Gamma \curvearrowright \Omega$, we can ask for a $\Gamma$-semilattice $E$, together with a $\Gamma$-equivariant isomorphism $\Zz[E\reg] \cong C_0(\Omega, \Zz)$. It is easy to see that such a system $\Gamma \curvearrowright E$ exists for $\Gamma \curvearrowright \Omega$ if and only if $\Omega$ admits a $\Gamma$-invariant regular basis in the sense of \cite[Definition~2.9]{C-E-L2}, i.e., a regular basis for $\Omega$ which is invariant under the $\Gamma$-action. In general, this does not need to be the case, as was remarked in \cite[Proposition~3.18]{C-E-L2}. However, an analogous argument as in \cite[Remark~3.22]{C-E-L2}, for integral semigroup rings instead of C*-algebras of semilattices, shows that given an arbitrary totally disconnected dynamical system $\Gamma \curvearrowright \Omega$, we can always find $\Gamma$-semilattices $E$, $E_1$, $E_2$, ..., and a $\Gamma$-equivariant long exact sequence
\bgl
\label{ind-res}
  \dotso \to \Zz[E_2\reg] \to \Zz[E_1\reg] \to \Zz[E\reg] \to C_0(\Omega, \Zz) \to 0.
\egl
We call such a long exact sequence an algebraic independent resolution of $\Gamma \curvearrowright C_0(\Omega, \Zz)$. Of course, the requirement that the sequence is $\Gamma$-equivariant is crucial here. Moreover, we define the length of such an algebraic independent resolution to be the smallest integer $n \geq 0$ with $E_{n+1} = \gekl{0}$, or equivalently, $\Zz[E_{n+1}\reg] = \gekl{0}$. If no such integer exists, then we set the length to be $\infty$.

As remarked at the beginning, we can also talk about independent resolutions of $\Gamma \curvearrowright D$ where $D$ is a commutative $\Zz$-algebra generated by idempotents with a $\Gamma$-action.

\bremark
Given a totally disconnected dynamical system $\Gamma \curvearrowright \Omega$, the minimal length of independent resolutions of $\Gamma \curvearrowright C_0(\Omega, \Zz)$ is an invariant of our dynamical system. The invariant is zero if and only if $\Omega$ admits a $\Gamma$-invariant regular basis in the sense of \cite[Definition~2.9]{C-E-L2}. This happens for instance for full shifts (see \cite[Example~3.1]{C-E-L2} for the case of two symbols, the general case of finitely many symbols is analogous). For the dynamical systems of \cite[Example~3.20]{C-E-L2}, this invariant is one. The reason is that it cannot be zero as explained in \cite[Example~3.20]{C-E-L2}, and the Toeplitz extension gives rise to an algebraic independent resolution of length one. It might be interesting to study this invariant from the point of view of dynamical systems. An obvious question would be what this invariant measures, maybe some sort of complexity? Another natural question would be whether all values in $\gekl{0,1,2,\dotsc} \cup \gekl{\infty}$ are actually possible. Also, it would be helpful to develop methods to compute this invariant, or at least to decide whether it is finite or infinite. For example, a concrete question would be: Is this invariant finite for shifts of finite type?
\eremark

\section{Existence of independent resolutions}

We have already argued that in principle, independent resolutions always exist. Now, we want to make this more precise.

Let $E$ be a semilattice. We introduce the following notation: Given finitely many $e_j \in E\reg$, $j \in J$, we form the idempotent $\bigvee_{j \in J} e_j \in \Zz[E\reg]$ given by the following sum:
$$
  \bigvee_{j \in J} e_j = \sum_{\varnothing \neq J' \subseteq J} (-1)^{\abs{J'}} \prod_{j \in J'} e_j.
$$
Note that $\bigvee_{j \in J} e_j$ is the smallest idempotent in $\Zz[E\reg]$ which dominates all the $e_j$. This is easy to see if $\abs{J} = 2$, since $e_1 \vee e_2 = e_1 + e_2 - e_1 e_2$, and the general case follows by an easy induction argument. This is the motivation for the formula above.

Let $D$ be a $\Zz$-algebra generated by commuting idempotents, and $\varphi$: $\Zz[E\reg] \to D$ a homomorphism of $\Zz$-algebras. Set
\bgloz
E_{\varphi} \defeq \menge{e - \bigvee_{j \in J} e_j \in \Zz[E\reg]}{J \ {\rm finite} \ {\rm set}, \, e, e_j \in E\reg, \, e_j \lneq e, \, \varphi(e) = \bigvee_{j \in J} \varphi(e_j)}.
\egloz

\blemma
{$\,$}
\begin{itemize}
\item[(i)] $E_{\varphi}$ is an ideal of $E$, i.e., $E \cdot E_{\varphi} \subseteq E_{\varphi}$,
\item[(ii)] $E_{\varphi}$ is multiplicatively closed,
\item[(iii)] $\ker \varphi = \lspan(E_{\varphi})$.
\end{itemize}
\elemma
By $\lspan$, we always mean $\Zz$-$\lspan$.
\bproof
(i) is verified by a straightforward computation. For (ii), take $e - \bigvee_{j \in J_e} e_j$ and $f - \bigvee_{j \in J_f} f_j$ in $E_{\varphi}$, and compute
\bglnoz
  && (e - \bigvee_{i \in J_e} e_i) (f - \bigvee_{j \in J_f} f_j) = ef - e \cdot (\bigvee_{j \in J_f} f_j) - (\bigvee_{i \in J_e} e_i) \cdot f + (\bigvee_{i \in J_e} e_i)(\bigvee_{j \in J_f} f_j) \\
  &=& ef - (\bigvee_{j \in J_f} e f_j) \vee (\bigvee_{i \in J_e} e_i f) \in E_{\varphi}.
\eglnoz
For (iii), let $I \defeq \lspan(E_{\varphi})$. $I$ is an ideal in $\Zz[E\reg]$. We have to show that $\varphi$ induces an injective homomorphism $\Zz[E\reg] / I \to D$. Let $\cF$ be the family of finite subsets $F$ of $E\reg$ which have the property that $F \cup \gekl{0}$ is multiplicatively closed. If we set $I_F \defeq I \cap \lspan(F)$ for $F \in \cF$, then we obviously have $\Zz[E\reg] / I \cong \ilim_{F \in \cF} \lspan(F) / I_F$. Here we order $\cF$ by inclusion. Therefore, it suffices to show that $\varphi$ induces an injective homomorphism on $\lspan(F) / I_F$ for each $F \in \cF$. For every $F \in \cF$, we get by orthogonalization a $\Zz$-basis for $\lspan(F)$ of pairwise orthogonal elements of the form $e - \bigvee_{i \in J} e_i$ (where $e_i \lneq e$). By construction of $I$, whenever $\varphi(e - \bigvee_{i \in J} e_i) = 0$, then $e - \bigvee_{i \in J} e_i$ lies in $I$, hence in $I_F$.
\eproof

\bcor
$\Zz[E_{\varphi}\reg] \lori \Zz[E\reg] \overset{\varphi}{\lori} D$ is exact, where the first homomorphism is induced by the canonical inclusion $E_{\varphi} \into \Zz[E\reg]$ (by universal property of $\Zz[E_{\varphi}\reg]$).
\ecor

\bremark
If a group $\Gamma$ acts on $E$ and $D$ such that $\varphi$ is $\Gamma$-equivariant, then the induced $\Gamma$-action on $\Zz[E\reg]$ restricts to a $\Gamma$-action on $E_{\varphi}\reg$. 
\eremark

The following observation will be useful.
\blemma
\label{fix-maximal-idem}
Let $E$ be a $\Gamma$-semilattice, where the $\Gamma$-action is given by $\Gamma \times E \to E, \, (g,e) \ma \tau_g(e)$. Let $e \in E\reg$ and $\gekl{e_j}_{j \in J}$ be a finite subset of $E\reg$ with $e_j \lneq e$ for all $j \in J$. For $g \in \Gamma$, $\tau_g(e - \bigvee_{j \in J} e_j) = e - \bigvee_{j \in J} e_j$ implies $\tau_g(e) = e$.
\elemma
\bproof
Since $E\reg$ is a $\Zz$-basis for $\Zz[E\reg]$, every element $x$ of $\Zz[E\reg]$ can be written in a unique way as $\sum_{\ve \in E\reg} n_{\ve} \ve$, where $n_{\ve}$ are integers, and only finitely many of the $n_{\ve}$s are not zero. We refer to this as the reduced form of $x$. In particular, for $x = \bigvee_{j \in J} e_j$, we can find a finite subset $F \subseteq E\reg$ and integer coefficients $n_{\ve}$ such that $\bigvee_{j \in J} e_j = \sum_{\ve \in F} n_{\ve} \ve$. The assumption $e_j \lneq e$ for all $j \in J$ implies $\ve \lneq e$ for all $\ve \in F$ since $F \subseteq \menge{\prod_{j \in J'} e_j}{\varnothing \neq J' \subseteq J}$ by the above formula for $\bigvee_{j \in J} e_j$. Hence $e - \sum_{\ve \in F} n_{\ve} \ve$ is the reduced form of $e - \bigvee_{j \in J} e_j$, and we get $\tau_g(e - \bigvee_{j \in J} e_j) = \tau_g(e) - \sum_{\ve \in F} n_{\ve} \tau_g(\ve)$, with $\tau_g(\ve) \lneq \tau_g(e)$ for all $\ve \in F$ because the group action preserves the relation $\lneq$. Since $\tau_g(e - \bigvee_{j \in J} e_j) = e - \bigvee_{j \in J} e_j$, we obtain by uniqueness of the reduced form that $\gekl{e} \cup F = \gekl{\tau_g(e)} \cup \tau_g(F)$. But $e$ is the unique element of $\gekl{e} \cup F$ which dominates all the elements in $\gekl{e} \cup F$, and similarly, $\tau_g(e)$ is the unique element of $\gekl{\tau_g(e)} \cup \tau_g(F)$ which dominates all the elements in $\gekl{\tau_g(e)} \cup \tau_g(F)$. We deduce that $e = \tau_g(e)$, as desired.
\eproof

\bprop
Let $D$ be a $\Zz$-algebra generated by commuting idempotents as above, and assume that a group $\Gamma$ acts on $D$. Then
\begin{itemize}
\item[(i)] There always exists a $\Gamma$-semilattice $E$ with a $\Gamma$-equivariant, surjective homomorphism $\varphi$: $\Zz[E\reg] \onto D$.
\item[(ii)] Choose $E$ and $\varphi$ with the properties in (i). If we set $E_0 \defeq E$, $\varphi_0 \defeq \varphi$, and recursively define $E_{k+1} \defeq E_{\varphi_k}$, $\varphi_{k+1}$ as the homomorphism $\Zz[E_{k+1}\reg] \to \Zz[E_k\reg]$ induced by $E_{k+1} \into E_k$, and if we equip all these semilattices with the canonical $\Gamma$-actions, then
\bgloz
  \dotso \to \Zz[E_k\reg] \to \dotso \to \Zz[E_2\reg] \to \Zz[E_1\reg] \to \Zz[E\reg] \to D \to 0
\egloz
is an algebraic independent resolution of $D$.
\item[(iii)] If $\Gamma$ is torsionfree, then there always exists a $\Gamma$-semilattice $E$ as in (i) such that the corresponding algebraic independent resolution in (ii) has the property that $\Gamma$ acts freely on $E_k\reg$.
\end{itemize}
\eprop
\bproof
For (i), just choose a $\Gamma$-invariant semilattice $E$ of idempotents in $D$ which generates $D$. (ii) follows from the corollary. To prove (iii), by the remark, it suffices to find a $\Gamma$-semilattice $E$ with $\Gamma$ acting freely on $E\reg$ and a $\Gamma$-equivariant homomorphism $\Zz[E\reg] \onto D$. Let us first treat the special case $D = \Zz$ together with the trivial $\Gamma$-action. Consider $E_{\Gamma} \defeq \menge{\Gamma \setminus F}{\varnothing \neq F \subseteq \Gamma \ \text{finite}} \cup \gekl{\varnothing}$. Here $\Gamma \setminus F$ is the complement of $F$ in $\Gamma$. $E_{\Gamma}$ is a semilattice under intersection. $\Gamma$ acts freely on $E_{\Gamma}\reg$ because $\Gamma$ is torsionfree. For arbitrary $X$ and $Y$ in $E_{\Gamma}$ not equal to $\varnothing$, we have that the product of $X$ and $Y$, $X \cap Y$, is not $\varnothing$. Therefore, there exists a semilattice homomorphism $E_{\Gamma} \to \gekl{0,1}$ determined by $X \ma 1$ for $X \neq \varnothing$ and $\varnothing \ma 0$. This homomorphism induces the homomorphism $\phi$: $\Zz[E_{\Gamma}\reg] \onto \Zz = \Zz[\gekl{0,1}\reg]$ which we need. Let us now turn to the general case. For $D$ as given, because of (i), we can always find $E$ and $\varphi$ with the properties in (i). Take the semilattice $E_{\Gamma}$ from above, and form the semilattice $E_{\Gamma} \times_0 E \defeq (E_{\Gamma}\reg \times E\reg) \cup \gekl{0}$ with the obvious multiplication. $\Gamma$ acts on $E_{\Gamma} \times_0 E$ diagonally, and hence freely on $(E_{\Gamma} \times_0 E)\reg = E_{\Gamma}\reg \times E\reg$, because the $\Gamma$-action on $E_{\Gamma}\reg$ is free. Moreover, we have $\Zz[(E_{\Gamma} \times_0 E)\reg] \cong \Zz[E_{\Gamma}\reg] \otimes_{\Zz} \Zz[E\reg]$. Thus $\phi \otimes \varphi$: $\Zz[E_{\Gamma}\reg] \otimes_{\Zz} \Zz[E\reg] \onto \Zz \otimes_{\Zz} D \cong D$ gives the desired homomorphism.
\eproof

\section{Existence of finite length independent resolutions}
\label{exist-fl-ind-res}

As explained in the introduction, we are particularly interested in independent resolutions of finite length. Given a commutative $\Zz$-algebra $D$ generated by idempotents together with a group action $\Gamma \curvearrowright D$, the only reason why we cannot find a $\Gamma$-semilattice $E$ such that $D$ is $\Gamma$-equivariantly isomorphic to $\Zz[E\reg]$ is that the idempotents in $D$ satisfy non-trivial linear relations over $\Zz$ which mix in a complicated way with the group action. Our goal now is to give conditions on these relations which allow us to find finite length independent resolutions. This requires algebraic manipulations in semilattices. We will often use orthogonalization as a tool. Later on, in \S~\ref{groups}, and in \cite[\S~6]{L-N}, we show that the conditions we introduce here are satisfied in natural examples, so that our Theorem~\ref{maintheo} is really useful in applications.

We need a bit of notation first. Let $E$ be a semilattice. A finite cover for $e \in E\reg$ is a finite subset $\gekl{f_j}_{j \in J}$ of $E\reg$ ($J$ is a finite index set) with the property that
\begin{itemize}
\item $f_j \leq e$ for all $j \in J$,
\item for every $f \in E\reg$ with $f \leq e$, there exists $j \in J$ such that $f f_j \neq 0$.
\end{itemize}
The definition of finite covers for elements of semilattices was introduced in \cite{Exel} and plays an important role in the context of tight representations. For our purpose it has an algebraic motivation. Let $D$ be a $\Zz$-algebra generated by commuting idempotents and $\varphi:\Zz[E\reg]\to D$ a homomorphism of $\Zz$-algebras. Let $\gekl{f_j}_{j \in J}\subset E\reg$ be such that $f_j\leq e$ for each $j\in J$ and $\varphi(e)=\bigvee_{j\in J}\varphi(f_j)$. If $\gekl{f_j}_{j \in J}$ is not a finite cover for $e$, then there is some nonzero $f\leq e$ such that $ff_j=0$ for each $j\in J$ and it follows that $\varphi(f)=0$. Hence in this situation, $E$ contains elements that are meaningless for the purpose of studying $D$. The term ``finite cover'' is also motivated by set theory. If $X$ is a non-empty set and $X_1,\ldots,X_n$ are subsets of $X$, then $\bigcup_{j=1}^n X_j=X$ if and only if for every non-empty subset $Y$ of $X$ there is a $j$ such that $Y\cap X_j\neq\varnothing$.

Let $\gekl{f_j}_{j \in J}$ be a finite cover for $e \in E\reg$. Recall that $\bigvee_{j \in J} f_j$ is the smallest idempotent in $\Zz[E\reg]$ dominating all the $f_j$, $j \in J$. Let $\bigvee\gekl{f_j}_{j \in J} \defeq \bigvee_{j \in J} f_j$. Furthermore, given another element $d \in E\reg$, we write $d \cdot \gekl{f_j}_{j \in J} \defeq \menge{d f_j}{j \in J} \eqdef \gekl{f_j}_{j \in J} \cdot d$ and $(d \cdot \gekl{f_j}_{j \in J})\reg \defeq (d \cdot \gekl{f_j}_{j \in J}) \cap E\reg = (\gekl{f_j}_{j \in J} \cdot d) \cap E\reg \eqdef (\gekl{f_j}_{j \in J} \cdot d)\reg$. Moreover, every $x \in \Zz[E\reg]$ can be written in a unique way as $\sum_{\ve \in E\reg} n_{\ve} \ve$, with finitely many $n_{\ve}$ not equal to zero. Let $E(x) \defeq \menge{\ve \in E\reg}{n_{\ve} \neq 0}$ be the support of $x$.

Our goal now is to produce an (explicit) algebraic independent resolution in the following situation: Let $E$ be a $\Gamma$-semilattice, where the $\Gamma$-action on $E$ is denoted by $\Gamma \times E \to E, (g,e) \ma \tau_g(e)$. Let us assume that we are given a collection of finite covers $\sR$ for $E$, i.e., for every $e \in E\reg$ a set $\sR(e)$ of finite covers for $e$, such that the following hold:
\begin{enumerate}
\item[(i)] For $d$, $e$ in $E\reg$ with $de \neq 0$ and $\cR \in \sR(e)$, either $de \in (d \cdot \cR)\reg$ or $(d \cdot \cR)\reg \in \sR(d e)$.
\item[(ii)] For $e \in E\reg$, pairwise distinct $\cR_1, \dotsc, \cR_r$ in $\sR(e)$ and $\ve_i \in E(\bigvee\cR_i)$ for $1 \leq i \leq r$, let $\ve \defeq \prod_{i=1}^r \ve_i$ and for $1 \leq j \leq r$, $\ve(\check{j}) \defeq \prod_{\substack{i=1 \\ i \neq j}}^r \ve_i$ (for $r = 1$, $\ve(\check{j}) = e$). Then we have for every $1 \leq j \leq r$: If $\ve(\check{j}) \neq 0$, then $\ve \lneq \ve(\check{j})$.
\item[(iii)] For every $g \in \Gamma$ and $e \in E\reg$, we have $\tau_g (\sR(e)) = \sR(\tau_g(e))$.
\end{enumerate}
We think of $\bigcup_{e \in E\reg} \sR(e)$ as a set of relations which we impose on the idempotents $e$ in $E\reg$, i.e., we want $e = \bigvee_{f \in \cR} f$ to hold for every $\cR \in \sR(e)$. Imposing such relations precisely corresponds to forming the quotient of $\Zz[E\reg]$ by the ideal $I$ generated by $e - \bigvee\cR$, $e \in E\reg$, $\cR \in \sR(e)$.

Roughly speaking, (i) and (iii) make sure that the relations we impose are compatible with the semigroup structure and the group action, whereas (ii) guarantees that the finite covers in $\sR$ interact in a very controlled way. We will see in the sequel why these conditions (i), (ii) and (iii) are useful. They also appear naturally, as will become clear in in \S~\ref{groups} and \cite[\S~6]{L-N}.

Obviously, (iii) ensures that the $\Gamma$-action induces a $\Gamma$-action on the quotient $\Zz[E\reg] / I$. We now set out to produce an (explicit) algebraic independent resolution for the dynamical system $\Gamma \curvearrowright \Zz[E\reg] / I$. The idea is to produce $\Gamma$-semilattices $E_1, E_2, \dotsc$ inductively and then to prove that they fit into a $\Gamma$-equivariant exact sequence. More precisely, we proceed as follows: Given a $\Gamma$-semilattice $E$ and a collection of finite covers $\sR$ for $E$ satisfying (i), (ii) and (iii), the first step is to construct a new $\Gamma$-semilattice $E(E,\sR)$ and a new collection of finite covers $\sR(E,\sR)$ for $E(E,\sR)$. The second step is to show that $E(E,\sR)$ and $\sR(E,\sR)$ satisfy (i), (ii) and (iii). In the third step, we construct the desired long exact sequence using an inductive argument.

\paragraph{\it First step}

We start with the construction of $E(E,\sR)$ and $\sR(E,\sR)$. Let $E$ be a $\Gamma$-semilattice, $\sR$ a collection of finite covers for $E$, and assume that $E$ and $\sR$ satisfy (i), (ii) and (iii).

\bdefin
Given $e \in E\reg$ and a non-empty finite subset $\sS_e$ of $\sR(e)$, set $e(\sS_e) \defeq \prod_{\cR \in \sR(e)} (e - \bigvee \cR) \in \Zz[E\reg]$. Define 
$$E(E,\sR) \defeq \menge{e(\sS_e)}{e \in E\reg, \, \varnothing \neq \sS_e \subseteq \sR(e) \ {\rm finite}} \cup \gekl{0}.$$

For $e \in E\reg$, a non-empty finite subset $\sS_e$ of $\sR(e)$ and $\ti{\cR} \in \sR(e) \setminus \sS_e$, set $\cR(\sS_e,\ti{\cR}) \defeq \gekl{e(\sS_e \cup \gekl{\ti{\cR}})} \cup \menge{f \cdot e(\sS_e)}{f \in \ti{\cR}, \, f \cdot e(\sS_e) \neq 0}$. Define
$$\sR(E,\sR)(e(\sS_e)) \defeq \menge{\cR(\sS_e,\ti{\cR})}{\ti{\cR} \in \sR(e) \setminus \sS_e}.$$
\edefin
The idea behind the construction of $E(E,\sR)$ and $\sR(E,\sR)$ is as follows: It is clear that $e - \bigvee \cR$ is in the kernel of the canonical projection $\Zz[E\reg] \onto \Zz[E\reg]/I$, and $E(E,\sR)$ is just the subsemigroup of $\Idem(\Zz[E\reg])$ generated by elements of the form $e - \bigvee \cR$. Moreover, $\sR(E,\sR)$ describe relations among the idempotents in $E(E,\sR)$ given by
\bgl
\label{evfe}
  e(\sS_e) = e(\sS_e \cup \gekl{\ti{\cR}}) \vee \bigvee_{f \in \ti{\cR}} f \cdot e(\sS_e)
\egl
for $\ti{\cR} \in \sR(e) \setminus \sS_e$.

Let us now show that $E(E,\sR)$ is a semilattice and that $\sR(E,\sR)$ is a collection of finite covers for $E(E,\sR)$.

\blemma
\label{1.L1}
$E(E,\sR)$ is a subsemigroup of $\Idem(\Zz[E\reg])$.
\elemma
\bproof
It is clear that $E(E,\sR) \subseteq \Idem(\Zz[E\reg])$. Let us show that $E(E,\sR)$ is multiplicatively closed. Let $d, e \in E\reg$, $\cQ \in \sR(d)$, $\cR \in \sR(e)$. All we have to show is that if $(d-\bigvee\cQ)(e-\bigvee\cR) \neq 0$, then
$$
  (d-\bigvee\cQ)(e-\bigvee\cR) = \prod_{\cS \in \sS_{de}} (de - \bigvee\cS)
$$
for some non-empty finite subset $\sS_{de}$ of $\sR(de)$.

Since $(d-\bigvee\cQ)(e-\bigvee\cR) \neq 0$, we know that $d \cdot (e- \bigvee \cR) \neq 0$. By condition (i), we must have $(d \cdot \cR)\reg = \cS$ for some $\cS \in \sR(de)$. This implies that $d \cdot \bigvee \cR = d \cdot \bigvee_{f \in \cR} f = \bigvee_{\ti{f} \in \cS} \ti{f} = \bigvee \cS$ and thus $d \cdot (e-\bigvee\cR) = de - \bigvee\cS$. Similarly, $(\cQ \cdot e)\reg = \cT$ for some $\cT \in \sR(de)$, and $(d-\bigvee\cQ) \cdot e = de - \bigvee\cT$. Thus, we conclude that
$$
  (d-\bigvee\cQ)(e-\bigvee\cR) = (d-\bigvee\cQ) \cdot e \cdot d \cdot (e-\bigvee\cR) = (de - \bigvee\cT) (de - \bigvee\cS).
$$
\eproof

\blemma
\label{1.L2}
Let $e \in E\reg$, $\sS_e$ be a non-empty finite subset of $\sR(e)$, and $\ti{\cR} \in \sR(e) \setminus \sS_e$. Then $\cR(\sS_e,\ti{\cR})$ is a finite cover for $e(\sS_e)$.
\elemma
\bproof
We first have to show that $\cR(\sS_e,\ti{\cR})$ is contained in $E(E,\sR)\reg$. It is clear that $e(\sS_e \cup \gekl{\ti{\cR}}) \in E(E,\sR)$. To show that it is not zero, we prove that for every $e \in E\reg$ and non-empty finite subset $\sS_e$ of $\sR(e)$, the idempotent $\prod_{\cR \in \sS_e}(e-\bigvee\cR) \in \Zz[E\reg]$ is not zero: If we expand the product $\prod_{\cR \in \sS_e}(e-\bigvee\cR)$, we obtain
$$
  \prod_{\cR \in \sS_e}(e-\bigvee\cR) = e - \sum_{\cR \in \sS_e} \bigvee\cR + \dotso.
$$
By definition, $\bigvee \cR$ is a linear combination of $\menge{\ve}{\ve \in E(\bigvee\cR)}$. Each of the remaining idempotents in the sum is below some $\ve \in E(\bigvee\cR)$, $\cR \in \sS_e$. Condition (ii) tells us that $\ve \lneq e$ for all $\ve \in E(\bigvee\cR)$. As the idempotents $e \in E\reg$ are linearly independent in $\Zz[E\reg]$, the product $\prod_{\cR \in \sS_e}(e-\bigvee\cR)$ cannot be zero. To prove that $\cR(\sS_e,\ti{\cR}) \subseteq E(E,\sR)\reg$, it remains to show that for every $f \in \ti{\cR}$, $f \cdot e(\sS_e) \in E(E,\sR)$. If $f \cdot (e-\bigvee \cR) = 0$ for some $\cR \in \sS_e$, then there is nothing to show. Otherwise, condition (i) tells us that $(f \cdot \cR)\reg = \cS_{\cR}$ for some $\cS_{\cR} \in \sR(f)$. Therefore, 
$$
  f \cdot \prod_{\cR \in \sS_e}(e-\bigvee\cR) = \prod_{\cR \in \sS_e} (f - \bigvee\cS_{\cR}) \in E(E,\sR).
$$
By construction, $e(\sS_e) = \bigvee_{f' \in \cR(\sS_e,\ti{\cR})} f'$ in $\Zz[E\reg]$ (see \eqref{evfe}). Therefore, $\cR(\sS_e,\ti{\cR})$ is a finite cover for $e(\sS_e)$.
\eproof

\paragraph{\it Second step}

Let us prove (i), (ii) and (iii) for $E(E,\sR)$ and $\sR(E,\sR)$. We need
\blemma
\label{2.L1}
Let $d, \, e \in E\reg$, $\sS_e$ be a non-empty finite subset of $\sR(e)$, $\ti{\cR} \in \sR(e) \setminus \sS_e$ and $\cQ \in \sR(d)$. If $(d-\bigvee\cQ)(e-\bigvee\cR) \neq 0$ for all $\cR \in \sS_e$, then $(d-\bigvee\cQ) \cdot e(\sS_e) \in (d-\bigvee\cQ) \cdot \cR(\sS_e,\ti{\cR}) $ or $((d-\bigvee\cQ) \cdot \cR(\sS_e,\ti{\cR}))\reg \in \sR(E,\sR)((d-\bigvee\cQ) \cdot e(\sS_e))$.
\elemma
\bproof
First assume that $(d-\bigvee\cQ)(e-\bigvee\ti{\cR}) = 0$. It follows from $(d-\bigvee\cQ)(e-\bigvee\cR) \neq 0$ for all $\cR \in \sS_e$ and our computations in Lemma~\ref{1.L1} and Lemma~\ref{1.L2} that we must have $d \cdot (e-\bigvee\ti{\cR}) = 0$. Since the idempotents in $E\reg$ are independent, we deduce that there must exist $f \in \ti{\cR}$ such that $de = df$. Therefore, $(d-\bigvee\cQ) \cdot e(\sS_e) = (d-\bigvee\cQ) \cdot f \cdot e(\sS_e)$ lies in $(d-\bigvee\cQ) \cdot \cR(\sS_e,\ti{\cR})$.

Secondly, assume that $(d-\bigvee\cQ)(e-\bigvee\cR) \neq 0$ for all $\cR \in \sS_e \cup \gekl{\ti{\cR}}$. As in the proof of Lemma~\ref{1.L1}, it follows that we can find $\cS, \ti{\cT}$ and $\cT_{\cR}$ in $\sR_{de}$ for every $\cR \in \sS_e$ such that $(\cQ \cdot e)\reg = \cS, \ (d \cdot \cR)\reg = \cT_{\cR} \fa \cR \in \sS_e, \ (d \cdot \ti{\cR})\reg = \ti{\cT}$. Thus $(d-\bigvee\cQ) \cdot e = de - \bigvee\cS$, $d \cdot (e-\bigvee\cR) = de - \bigvee\cT_{\cR}$ and $d \cdot (e-\bigvee\ti{\cR}) = de - \bigvee\ti{\cT}$.

If $\ti{\cT}$ lies in $\gekl{\cS} \cup \menge{\cT_{\cR}}{\cR \in \sS_e}$, then
\bglnoz
  (d-\bigvee\cQ) \cdot e(\sS_e) = (d-\bigvee\cQ) \cdot e(\sS_e \cup \gekl{\ti{\cR}}) \in (d-\bigvee\cQ) \cdot \cR(\sS_e,\ti{\cR}).
\eglnoz

If $\ti{\cT} \notin \gekl{\cS} \cup \menge{\cT_{\cR}}{\cR \in \sS_e}$, then
\bgloz
  (d-\bigvee\cQ) \cdot e(\sS_e) = (de)(\menge{\cT_{\cR}}{\cR \in \sS_e} \cup \gekl{\cS})
\egloz
and
\bgloz
  (d-\bigvee \cQ) \cdot e(\sS_e \cup \gekl{\ti{\cR}}) = (de)(\menge{\cT_{\cR}}{\cR \in \sS_e} \cup \gekl{\cS,\ti{\cT}}).
\egloz
As $(d \cdot \ti{\cR})\reg = \ti{\cT}$, we have
\bglnoz
  && \menge{(d-\bigvee\cQ) \cdot f \cdot e(\sS_e)}{f \in \ti{\cR}}\reg \\
  &=& \menge{d \cdot f \cdot (de)(\menge{\cT_{\cR}}{\cR \in \sS_e} \cup \gekl{\cS})}{f \in \ti{\cR}}\reg \\
  &=& \menge{\ti{f} \cdot (de)(\menge{\cT_{\cR}}{\cR \in \sS_e} \cup \gekl{\cS})}{\ti{f} \in \ti{\cT}},
\eglnoz
and thus
$$((d - \bigvee \cQ) \cdot \cR(\sS_e,\ti{\cR})\reg = \cR(\menge{\cT_{\cR}}{\cR \in \sS_e} \cup \gekl{\cS}, \ti{\cT}) \in \sR(E,\sR)((d-\bigvee\cQ) \cdot e(\sS_e)).$$
\eproof

\blemma
\label{2.L2}
$E(E,\sR)$ and $\sR(E,\sR)$ satisfy condition (i).
\elemma
\bproof
Let $e(\sR)$ and $\ti{\cR}$ be as in Lemma~\ref{2.L1}, and let $d' = d(\gekl{\cQ_1, \dotsc, \cQ_m})$ for $d \in E\reg$, $\cQ_1, \dotsc, \cQ_m \in \sR(d)$. We show inductively on $m$ that $d' e(\sS_e) \in d' \cdot \cR(\sS_e,\ti{\cR})$ or that $(d' \cdot \cR(\sS_e,\ti{\cR}))\reg \in \sR(E,\sR)(d' e(\sS_e))$ (under the assumption that $d' e(\sS_e) \neq 0$). For $m=1$ this is precisely the content of Lemma~\ref{2.L1}. Assume that we have proven our assertion for $d'' = d(\gekl{\cQ_1, \dotsc, \cQ_{m-1}})$. If $d'' e(\sS_e) \in d'' \cdot \cR(\sS_e,\ti{\cR})$, then $d' e(\sS_e) = (d-\bigvee\cQ_m) \cdot d'' \cdot e(\sS_e) \in (d-\bigvee\cQ_m) \cdot (d'' \cdot \cR(\sS_e,\ti{\cR})) = d' \cdot \cR(\sS_e,\ti{\cR})$. If  $(d'' \cdot \cR(\sS_e,\ti{\cR}))\reg \in \sR(E,\sR)(d'' e(\sS_e))$, then $(d' \cdot \cR(\sS_e,\ti{\cR}))\reg = ((d-\bigvee\cQ_m) \cdot (d'' \cdot \cR(\sS_e,\ti{\cR}))\reg)\reg$ either contains $(d-\bigvee\cQ_m) \cdot d'' \cdot e(\sS_e) = d' e(\sS_e)$ or is an element of $\sR(E,\sR)((d-\bigvee\cQ_m) \cdot d'' \cdot e(\sS_e)) = \sR(E,\sR)(d' e(\sS_e))$ by Lemma~\ref{2.L1}.
\eproof

Let us turn to (ii).
\blemma
\label{2.L3}
Given $e \in E\reg$ and $\cR_1, \dotsc, \cR_r \in \sR(e)$ pairwise distinct, choose for every $1 \leq i \leq r$ an idempotent $\ve_i \in E(\bigvee\cR_i) \cup \gekl{e-\bigvee\cR_i}$. Let $\ve \defeq \prod_{i=1}^r \ve_i$ and for $1 \leq j \leq r$, $\ve(\check{j}) \defeq \prod_{\substack{i=1 \\ i \neq j}}^r \ve_i$ (for $r = 1$, $\ve(\check{j}) = e$). Let $1 \leq j \leq r$. If $\ve(\check{j}) \neq 0$, then $\ve \lneq \ve(\check{j})$.
\elemma
\bproof
Assume $\ve(\check{j}) \neq 0$. The reduced form of $\ve(\check{j})$ is given by a linear combination of idempotents in $E\reg$. Among these non-zero idempotents, there will be a biggest one which appears with non-zero coefficient in this linear combination. All the other idempotents in our linear combination will be strictly smaller because of condition (ii). Let $d \in E\reg$ be this biggest idempotent. It is clear that $d$ is the product of those $\ve_i$s which were chosen from $E(\bigvee\cR_i)$, $i \neq j$.

First of all, let us assume that $\ve_j \in E(\bigvee\cR_j)$. The reduced form of $\ve = \prod_{i=1}^r \ve_i = \ve_j \cdot \ve(\check{j})$ is a linear combination of the same form as for $\ve(\check{j})$, but this time the biggest idempotent is $\ve_j \cdot d \neq 0$. By condition (ii), we must have $d \gneq \ve_j \cdot d$. Thus $\ve \lneq \ve(\check{j})$.

Next, assume that $\ve_j = e - \bigvee\cR_j$. If $\ve = \ve(\check{j})$, then
$$
  \ve(\check{j}) - \bigvee\cR_j \cdot \ve(\check{j}) = \ve(\check{j}).
$$
Thus $\bigvee\cR_j \cdot \ve(\check{j}) = 0$. Therefore, we must have
\bgl
\label{contra}
  \delta \cdot \ve(\check{j}) = 0 \fa \delta \in E(\bigvee\cR_j).
\egl
Recall that $d \neq 0$ is the biggest idempotent in the reduced form of $\ve(\check{j})$. We have $0 \neq d \leq e$. Since $\cR_j$ is a cover for $e$, there must be $\delta \in E(\bigvee\cR_j)$ such that $\delta \cdot d \neq 0$. As we have already seen, $\delta \cdot d$ is the biggest idempotent in the reduced form of $\delta \cdot \ve(\check{j})$. As $\delta \cdot d \neq 0$ and all the remaining idempotents in the expansion of $\delta \cdot \ve(\check{j})$ are strictly smaller than $\delta \cdot d$ by condition (ii), we conclude that $\delta \cdot \ve(\check{j}) \neq 0$. But this contradicts \eqref{contra}.
\eproof

\blemma
\label{2.L4}
$E(E,\sR)$ and $\sR(E,\sR)$ satisfy condition (ii).
\elemma
\bproof
This is an immediate consequence of Lemma~\ref{2.L3}. The reason is that every $\ve' \in E(E,\sR)(\bigvee_{i=1}^r \cR'_i)$, where $\cR'_i \in \sR(E,\sR)(e(\sS_e))$ are pairwise distinct, is of the form $\prod_i \ve_i$ as in Lemma~\ref{2.L3}.
\eproof

Let us discuss (iii). First of all, it is obvious that the $\Gamma$-action $\tau$ on $E$ induces a $\Gamma$-action on $E(E,\sR)$. More precisely, $\tau$ induces a $\Gamma$-action on $\Zz[E\reg]$ (again denoted by $\tau$), and we have because of $\tau_g(\sR_e) = \sR_{\tau_g(e)}$ (condition (iii)):
$$
  \tau_g(\prod_{\cR \in \sS_e}(e-\bigvee\cR)) = \prod_{\cR \in \sS_e}(\tau_g(e) - \tau_g(\bigvee\cR)) = \prod_{\cS \in \tau_g(\sS_e)}(\tau_g(e) - \bigvee\cS),
$$
so that $\tau$ restricts to a $\Gamma$-action $\tau'$ on $E(E,\sR) \subseteq \Idem(\Zz[E\reg])$.

\blemma
We have $\tau'_g(\sR(E,\sR)(e')) = \sR(E,\sR)(\tau'_g(e'))$ for all $e' \in E(E,\sR)\reg$, $g \in \Gamma$. In other words, the $\Gamma$-semilattice $E(E,\sR)$ and $\sR(E,\sR)$ satisfy (iii).
\elemma
\bproof
Let $e \in E\reg$ and $e' = e(\sS_e)$ for a non-empty finite subset $\sS_e$ of $\sR(e)$. Let $\ti{\cR} \in \sR(e) \setminus \sS_e$. Then (iii) implies that $\tau_g(\sS_e)$ is a a non-empty finite subset of $\sR(\tau_g(e))$, and that $\tau_g(\ti{\cR}) \in \sR(\tau_g(e)) \setminus \tau_g(\sS_e)$. Hence, by construction of $\cR(\sS_e,\ti{\cR})$, we have $\tau'_g(\cR(\sS_e,\ti{\cR}) = \cR(\tau_g(\sS_e),\tau_g(\ti{\cR})) \in \sR(E,\sR)(\tau_g(e)(\tau_g(\sS_e)))$.
\eproof

To summarize, we have proven
\bprop
\label{induction_i-iii}
$E(E,\sR)$ is a $\Gamma$-semilattice (with respect to) $\tau'$) and $\sR(E,\sR)$ is a collection of finite covers for $E(E,\sR)$ satisfying (i), (ii) and (iii).
\eprop

\paragraph{\it Third step}

To construct the desired long exact sequence, the following observation plays a crucial role. Let
\bglnoz
  && I = \spkl{\menge{e - \bigvee\cR}{e \in E\reg, \ \cR \in \sR(e)}}_{\Zz} \triangleleft \Zz[E\reg] \\
  && I' = \spkl{\menge{e' - \bigvee\cR'}{e' \in E(E,\sR)\reg, \ \cR' \in \sR(E,\sR)(e')}}_{\Zz} \triangleleft \Zz[E(E,\sR)\reg].
\eglnoz
By universal property of $\Zz[E(E,\sR)\reg]$, the inclusion $E(E,\sR) \into \Idem(\Zz[E\reg])$ induces a homomorphism
$$
  \pi: \ \Zz[E(E,\sR)\reg] \to \Zz[E\reg], \ E(E,\sR)\reg \ni e' = e(\sS_e) \ma \prod_{\cR \in \sS_e} (e-\bigvee\cR) \in \Zz[E\reg].
$$

\bprop
\label{induction_exact}
$ \ker \pi = I'$.
\eprop
\bproof
First of all, for $w \in \Zz[E(E,\sR)\reg]$, we set $E(w) = \bigcup_{e' \in E(E,\sR)(w)} E(e')$. We denote by $\bar{E}(w)$ the smallest multiplicatively closed subset of $E$ containing $E(w) \cup \gekl{0}$. Obviously, $\abs{\bar{E}(w)} < \infty$. For $e \in E(w)$, set $N_e(w) = \sum_{e' \in E(E,\sR)(w), \ e \in E(e')} \abs{n_{e'}}$. With these notations, we start the proof of our proposition.

Take $x \in \Zz[E(E,\sR)\reg]$ with $\pi (x) = 0$. We have to show that $x \in I'$. The idea is to proceed inductively on $\abs{\bar{E}(x)}$. Clearly, if $\abs{\bar{E}(x)} = 1$, then $\bar{E}(x) = \gekl{0}$, thus $E(x) \subseteq \gekl{0}$, and so $x=0$. Our goal is to find $z \in I'$ such that $\bar{E}(x-z) \subsetneq \bar{E}(x)$. Then by induction hypothesis, $x-z$ lies in $I'$, so $x$ lies in $I'$, and we are done.

Choose $e \in \bar{E}(x)$ maximal. $e$ must lie in $E(x)$ as every idempotent in $\bar{E}(x)$ is dominated by a idempotent in $E(x)$. Moreover, if $e \in E(e')$ for some $e' \in E(E,\sR)(x)$, then $e'$ must be of the form $e(\sS_e)$ since $e$ is maximal. We want to find $y \in I'$ such that $\bar{E}(x-y) \subseteq \bar{E}(x)$ and $N_e(x-y) < N_e(x)$. Since $\pi(x) = 0$ in $\Zz[E\reg]$, we must have $\sum_{e' \in E(E,\sR)(x), \ e \in E(e')} n_{e'} (x) = 0$. So there exist distinct $e'_1, e'_2 \in E(E,\sR)(x)$ with coefficients $n_i \defeq n_{e'_i}(x)$ such that $e \in E(e'_i)$ and $n_1 > 0$, $n_2 < 0$. We know that $e'_i$ is of the form $e'_i = e(\sS_i)$ with $\varnothing \neq \sS_i \subseteq \sR(e)$ finite ($i=1,2$), $\sS_1 \neq \sS_2$. Set $\sS = \sS_1 \cup \sS_2$ and write $\sS \setminus \sS_1 = \gekl{\cR_1^{(j)}}$, $\sS \setminus \sS_2 = \gekl{\cR_2^{(k)}}$. We obviously have the following identities modulo $I'$:
\bglnoz
  && e'_1 \\
  &\equiv& e(\sS_1 \cup \gekl{\cR_1^{(1)}}) + \bigvee_{f_1^{(1)} \in \cR_1^{(1)}} f_1^{(1)} \cdot e(\sS_1) \\
  &\equiv& e(\sS_1 \cup \gekl{\cR_1^{(1)}, \, \cR_1^{(2)}}) 
  + \bigvee_{f_1^{(2)} \in \cR_1^{(2)}} f_1^{(2)} \cdot e(\sS_1 \cup \gekl{\cR_1^{(1)}}) 
  + \bigvee_{f_1^{(1)} \in \cR_1^{(1)}} f_1^{(1)} \cdot e(\sS_1) \\
  &\equiv& \dots \\
  &\equiv& e(\sS) 
  + \underbrace{\sum_j \bigvee_{f_1^{(j)} \in \cR_1^{(j)}} f_1^{(j)} \cdot e(\sS_1 \cup \gekl{\cR_1^{(1)}, \dotsc, \cR_1^{(j-1)}})}_{\Delta_1}.
\eglnoz
Similarly, 
\bgloz
  e'_2 \equiv e(\sS)
  + \underbrace{\sum_k \bigvee_{f_2^{(k)} \in \cR_2^{(k)}} f_2^{(k)} \cdot e(\sS_1 \cup \gekl{\cR_2^{(1)}, \dotsc, \cR_2^{(k-1)}})}_{\Delta_2} \mod I'.
\egloz
So $y \defeq e'_1 - e'_2 - \Delta_1 + \Delta_2$ lies in $I'$. We claim that $\bar{E}(x-y) \subseteq \bar{E}(x)$ and $N_e(x-y) < N_e(x)$. To prove the first claim, note that since $e'_1, e'_2 \in E(E,\sR)(x)$, we must have $E(e'_1) \cup E(e'_2) \subseteq \bar{E}(x)$. Thus for every $\cR \in \sS$ and $\ve \in E(\bigvee\cR)$, $\ve$ lies in $\bar{E}(x)$. This makes use of condition (ii). By construction, $\bar{E}(x)$ is multiplicatively closed. Thus for every $\ti{\sS} \subseteq \sS$ and every choice of $\ve_{\cR} \in E(e(\cR))$, $\cR \in \ti{\sS}$, the product $\prod_{\cR \in \ti{\sS}} \ve_{\cR}$ lies in $\bar{E}(x)$. This shows that $\bar{E}(y) \subseteq \bar{E}(x)$. Hence $\bar{E}(x-y) \subseteq \bar{E}(x)$. For the second claim, note that $N_e(x-y) \leq N_e(x - e'_1 + e'_2) + N_e(\Delta_1 - \Delta_2) = N_e(x)-2 + 0 < N_e(x)$. So we have found $y \in I'$ with $\bar{E}(x-y) \subseteq \bar{E}(x)$ and $N_e(x-y) < N_e(x)$. After finitely many steps, we arrive at $y_1, \dotsc, y_m \in I'$ with $\bar{E}(x - y_1 - y_2 - \dotso - y_m) \subseteq \bar{E}(x)$ and $N_e(x - y_1 - y_2 - \dotso - y_m) = 0$. Thus with $z \defeq y_1 + \dotso + y_m \in I'$, we must have $\bar{E}(x-z) \subseteq \bar{E}(x)$ but $e \notin \bar{E}(x-z)$, hence $\bar{E}(x-z) \subsetneq \bar{E}(x)$.

Therefore we have found an element $z \in I'$ with the desired properties. This finishes the proof of the proposition.
\eproof

Finally, we are ready to prove our main result.
\btheo
\label{maintheo}
Suppose that $E$ is a $\Gamma$-semilattice and $\sR$ is a collection of finite covers for $E$ satisfying (i), (ii) and (iii). Let $I = \spkl{\menge{e - \bigvee \cR}{e \in E\reg, \ \cR \in \sR(e)}}_{\Zz} \triangleleft \Zz[E\reg]$. We define iteratively $E_0 \defeq E$, $\sR_0 \defeq \sR$ and for all $k \in \Nz$, $E_{k+1} \defeq E(E_k, \sR_k)$ and $\sR_{k+1} \defeq \sR(E_k, \sR_k)$ (using the notation from above). Then the canonical projection $\Zz[E\reg] \onto \Zz[E\reg] / I$ and the homomorphisms $\Zz[E_{k+1}\reg] \to \Zz[E_k\reg]$ induced by the inclusions $E_{k+1} \into \Idem(\Zz[E_k\reg])$ give rise to a $\Gamma$-equivariant long exact sequence
\bgl
\label{long-exact-seq}
  \dotso \to \Zz[E_2\reg] \to \Zz[E_1\reg] \to \Zz[E\reg] \to \Zz[E\reg] / I \to 0.
\egl
\etheo
\bproof
Proposition~\ref{induction_i-iii} tells us that we can really define $E_k$ and $\sR_k$ iteratively. To prove exactness of \eqref{long-exact-seq}, set $I_{k} = \spkl{\menge{e - \bigvee\cR}{e \in E_k\reg, \ \cR \in \sR_k(e)}}_{\Zz} \triangleleft \Zz[E_k\reg]$ for all $k \in \Nz$. The kernel of $\Zz[E\reg] \to \Zz[E\reg] / I$ is obviously given by $I_{0} = I$. Therefore, by induction on $k$, we obtain using Proposition~\ref{induction_exact} that the kernel of the canonical homomorphism $\Zz[E_{k+1}\reg] \to \Zz[E_k\reg]$ is given by $I_{k+1}$, and $I_{k}$ is the image of $\Zz[E_{k+1}\reg] \to \Zz[E_k\reg]$ by construction. Hence it follows that \eqref{long-exact-seq} is exact. 
\eproof

The following observations are immediate corollaries of our construction:
\bcor
\label{cor-finitelength}
In the situation of the theorem, if we have $\sup_{e \in E\reg} \abs{\sR(e)} < \infty$, then \eqref{long-exact-seq} is an algebraic independent resolution of $\Zz[E\reg] / I$ of length at most $\sup_{e \in E\reg} \abs{\sR(e)}$.
\ecor
\bproof
This follows immediately from the observation that $\abs{\sR(E,\sR)(e(\sS_e))} < \abs{\sR(e)}$ for every $e \in E\reg$ and every non-empty finite subset $\sS_e$ of $\sR(e)$.
\eproof

\bcor
\label{cor-stabilizer}
In the situation of the theorem, let $\Gamma^{(stab)}$ be the subgroup of $\Gamma$ generated by the stabilizer groups $\Gamma_e$, $e \in E\reg$. Then for every $n$ and $e' \in E_k\reg$, the stabilizer group $\Gamma_{e'}$ is contained in $\Gamma^{(stab)}$.
\ecor
\bproof
This is an immediate consequence of the following observation: Given $g \in \Gamma$, $e \in E\reg$ and a non-empty finite subset $\sS_e$ of $\sR(e)$, set $\Gamma_e \defeq \menge{g \in \Gamma}{\tau_g(e) = e}$ and $\Gamma_{e(\sS_e)} \defeq \menge{g \in \Gamma}{\tau'_g(e(\sS_e)) = e(\sS_e)}$. Then a similar argument as in the proof of Lemma~\ref{fix-maximal-idem} yields $\Gamma_{e(\sS_e)} \subseteq \Gamma_e$.
\eproof

\section{Computing homology in the case of free actions}
\label{compute-homology-free-action}

When we are in the situation that $\Gamma$ acts freely on the semilattice $E\reg$, the method of constructing independent resolutions presented in the previous section allows us to obtain a free resolution of $\Zz[E\reg] / I$ over $\Zz\Gamma$. This section is concerned with finding methods for explicitly computing the homology groups $H_*(\Gamma,\Zz[E\reg] / I)$. These homology groups will play a role in the computation of K-groups in \cite{L-N}, and may also be of independent interest. We first give a description of the chain complex $C$ that is derived directly from the construction of the independent resolution. We then introduce a chain complex $\wt{C}$ which is easier to use for homology computations, and which in some cases has the same homology as $C$.

Recall that in order to compute the group homology $H_*(\Gamma,A)$ of a $\Zz\Gamma$-module $A$, one finds a projective resolution of $A$ over $\Zz\Gamma$, i.e. a long exact sequence
\[
\cdots\to A_2\to A_1\to A\to 0
\]
where each $A_k$ is a projective $\Zz\Gamma$-module. Then one deletes $A$ from the resolution, and tensors the remaining sequence with the trivial $\Zz\Gamma$-module $\Zz$. The group homology is then the homology of the resulting chain complex (see for instance \cite{Bro}). Every free resolution is a projective resolution. Recall also that if $\Gamma$ acts freely on a set $X$, then there is a natural isomorphism $\oplus_X \Zz \cong \oplus_{\Gamma\setminus X}\Zz\Gamma$ of $\Zz\Gamma$-modules. This implies that $(\oplus_X \Zz) \otimes_{\Zz\Gamma}\Zz\cong \oplus_{\Gamma\setminus X}\Zz$.

Let $E$ be a fixed $\Gamma$-semilattice with a fixed system $\sR$ of covers satisfying (i)-(iii) of \S~\ref{exist-fl-ind-res}. Suppose also that $\Gamma$ acts freely on $E$. Then Corollary~\ref{cor-stabilizer} implies that $\Gamma$ acts freely on $E_k\reg$ for every $k \in \Nz$. So
\bgl\label{free-resolution--ZEI}
  \dotso \to \Zz[E_2\reg] \to \Zz[E_1\reg] \to \Zz[E\reg] \to \Zz[E\reg] / I \to 0.
\egl
is a free resolution of $\Zz[E\reg] / I$ over $\Zz\Gamma$. Deleting $\Zz[E\reg] / I$ from \eqref{free-resolution--ZEI} and tensoring with the trivial $\Zz\Gamma$-module we get a chain complex $C$ with $H_*(C)=H_*(\Gamma,\Zz[E\reg] / I)$. Computing these homology groups will be important in some applications. If we let $\Zz[\Gamma\setminus E_k\reg]$ be the free $\Zz$-module over the orbit space $\Gamma\setminus E_k\reg$, then $C$ can be viewed as the chain complex
$$
  C=\left(\dotso \to \Zz[\Gamma\setminus E_2\reg] \to \Zz[\Gamma\setminus E_1\reg] \to \Zz[\Gamma\setminus E\reg] \to 0\right)
$$
where the boundary maps are induced from the maps in \eqref{free-resolution--ZEI}. We will denote the boundary maps by $\partial_k:\Zz[\Gamma\setminus E_k\reg]\to\Zz[\Gamma\setminus E_{k-1}\reg]$. In order to compute homology, we need to describe these maps.

For $n,k\in\Nz$, $1\leq k\leq n$, let
$$
Q^n_k=\menge{(\mu_1,\ldots,\mu_k)\in (2^{\{1,\ldots,n\}}\setminus\{\varnothing\})^k}{\mu_i\cap \mu_j=\varnothing\mbox{ for }i\neq j}.
$$
So $Q^n_k$ is the set of $k$-tuples of mutually disjoint non-empty subsets of $\{1,\ldots,n\}$. Let $Q^n_0=\{\varnothing\}$. We will sometimes use a special notation for the elements of $Q^n_k$. Given $\mu\in Q^n_k$, write $\mu_i$ as a comma-separated list of its elements and write $\mu$ as a bar-separated list of the $\mu_i$s. For instance, given $\mu=(\{10,3\},\{4,5,8\},\{7\},\{1,2\})\in Q^{10}_4$, write $\mu=10,3|4,5,8|7|1,2$. Given $p\in\{1,\ldots n\}$ we say that $p\in\mu$ if $p\in\mu_i$ for some $i$. If $p\notin\mu$, $\mu|p$ stands for $(\mu_1,\ldots,\mu_k,\{p\})\in Q^n_{k+1}$. Shorten $\varnothing|p$ to $p$. Similarly, $\mu,p$ stands for $(\mu_1,\ldots,\mu_k\cup\{p\})\in Q^n_k$, and $\varnothing,p$ shortens to $p$. For $1\leq i\leq k$, let $\mu^i=(\mu_1,\ldots,\mu_{i-1},\mu_{i+1},\ldots,\mu_k)\in Q^n_{k-1}$. Given $\rho\subseteq\{1,\ldots,n\}$ such that $\rho\cap\mu_j=\varnothing$ for each $j$, let $\mu(i;\rho)=\mu_1|\cdots|\mu_{i-1}|\mu_i,\rho|\mu_{i+1}|\cdots|\mu_k$.

Let $0\leq k_1\leq k_2\leq n$, $\mu\in Q^n_{k_1}, \nu \in Q^n_{k_2}$. We construct a commutative partially defined product $\mu * \nu$ by the rule
$$
\mu * \nu = \begin{cases}
	\nu & k_1=0\\
	(\mu_1\cup\nu_1,\ldots,\mu_{k_1}\cup\nu_{k_1},\nu_{k_1+1},\ldots,\nu_{k_2}) & \mbox{when this lies in }Q^n_{k_2}\\
	\mbox{undefined} & \mbox{otherwise.}
	\end{cases}
$$
Note that the product is defined exactly when each $p\in\{1,\ldots n\}$ appears in at most one of the sets $\mu_1\cup\nu_1,\ldots,\mu_{k_1}\cup\nu_{k_1},\nu_{k_1+1},\ldots,\nu_{k_2}$. 

Let $E$ be a semilattice with a system $\sR$ of covers satisfying conditions (i)-(iii) of \S~\ref{exist-fl-ind-res}. Let $e\in E\reg$, and suppose $\sR(e)=\{\cR_1(e),\ldots, \cR_n(e)\}$. Given $\mu\in Q^n_k$, define $e(\mu)\in E_k$ iteratively using the following rules:
\begin{align*}
e(\varnothing)&=e\\
e(\mu|p)&=e(\mu)-\bigvee\cR_p(e(\mu))\in E_{k+1} & \mu\in Q^n_k, p\notin\mu\\
e(\mu|p_1,\ldots,p_m)&=\prod_{i=1}^m e(\mu|p_i) & \mu\in Q^n_k, p_i\notin\mu,1\leq i\leq m\\
\cR_p(e(\mu))&=\{e(\mu,p)\}\cup e(\mu)\cR_p(e(\mu^k)) & \mu\in Q^n_k, p\notin\mu.
\end{align*}
Note now that for $\mu\in Q^n_k$, $\sR_k (e(\mu)) = \menge{\cR_p(e(\mu))}{p\notin\mu}$. By the construction of $E_k$ we then have that if $|\sR(e)|$ is finite for each $e\in E\reg$,
$$
  E_k=\menge{e(\mu)}{e\in E\reg, n=|\sR(e)|, 1\leq k\leq n, \mu\in Q^n_k}.
$$
Let $k_1\geq k_2$, $\mu\in Q^n_{k_1},\nu\in Q^n_{k_2}$ and $\ve\in E(\bigvee\cR_j(e(\nu))$. We want to make sense of a product $e(\mu)\ve$. As in Lemma~\ref{1.L2} this immediately makes sense when $k_1=k_2+1$. We can then iteratively define
$$
e(\mu)\ve= \prod_{p\in\mu_{k_1}}\left(e(\mu^{k_1})\ve-\bigvee \cR_p(e(\mu^{k_1}))\ve \right),
$$
where $\cR_p(e(\mu^{k_1}))\ve=\{e(\mu^{k_1},p)\ve\}\cup (e(\mu^{k_1})\ve)(\cR_p(e((\mu^{k_1})^{k_1-1})\ve))$. It can be checked that $e(\mu)\ve$ is a well defined element of $E_{k_1}$. We also iteratively define
$$
e(\nu)e(\mu)=\prod_{p\in\mu_{k_1}}\left(e(\nu)e(\mu^{k_1})-\bigvee e(\nu)\cR_p(e(\mu^{k_1}))\right).
$$

\blemma\label{semilatticeproduct-higher-level}
Let $0\leq k_1,k_2\leq n$, $\mu\in Q^n_{k_1}, \nu \in Q^n_{k_2}$. We have
$$
e(\mu)e(\nu)=\begin{cases}
	e(\mu * \nu)&\mbox{if }\mu * \nu\mbox{ is defined}\\
	0&\mbox{otherwise}
\end{cases}
$$
Moreover, $e(\mu)\ve=0$ for any $\ve\in E(\bigvee\cR_j(e))$ with $j\in\mu$. If $j\notin\mu,\nu$, then $e(\mu)(e(\nu)\ve)=e(\mu * \nu)\ve$.
\elemma
\bproof
A critical step is to prove that given $1\leq i\leq n$, $i\leq k\leq n$, $\mu\in Q^n_k$ and $p,q\notin\mu$ we have
\begin{align}\label{eq:semilatticeproducthypothesis}
e(\mu_1|\cdots|\mu_{i-1}|p)e(\mu)&=e(\mu(i;p))\\
e(\mu_1|\cdots|\mu_{i-1}|p)\cR_q(e(\mu^k))&=\cR_q(e(\mu^k(i;p)))\mbox{ when }k>i\label{eq:semilatticeproducthypothesis2}
\end{align}
This is done by fixing $i$ and doing an induction proof with $k$ as the induction index. Suppose first that $k=i$. It follows directly from the definition that $e(\mu_1|\cdots|\mu_{i-1}|p)e(\mu)=e(\mu_1|\cdots|\mu_k,p)=e(\mu(i;p))$, so \eqref{eq:semilatticeproducthypothesis} holds. Then as $\mu^k(i;p)=\mu^k$,
\begin{align*}
&e(\mu_1|\cdots|\mu_{i-1}|p)\cR_q(e(\mu))=\\
&\{e(\mu_1|\cdots|\mu_{i-1}|p)e(\mu_1|\cdots|\mu_{i-1}|\mu_i,q)\}\cup e(\mu_1|\cdots|\mu_{i-1}|p)e(\mu)\cR_q(e(\mu^k))\\
&=\{e(\mu(i;p),q)\}\cup e(\mu(i;p)\cR_q(e(\mu^k(i;p)))=\cR_q(e(\mu(i;p))).
\end{align*}
So \eqref{eq:semilatticeproducthypothesis} holds. Suppose that now that the induction hypothesis holds for any $i\leq k\leq m$, $\mu\in Q^n_k$ and $p,q\notin\mu$ and show that this implies that it holds for any $k=m+1$, $\mu\in Q^n_k$ and $p,q\notin\mu$. We have
\begin{align*}
&e(\mu_1|\cdots|\mu_{i-1}|p)e(\mu)\\
&=\prod_{q\in\mu_k}\left(e(\mu_1|\cdots|\mu_{i-1}|p)e(\mu^k)-\bigvee e(\mu_1|\cdots|\mu_{i-1}|p)\cR_q(e(\mu^k))\right)\\
&= \prod_{q\in\mu_k}\left(e(\mu^k(i;p))-\bigvee\cR_{\mu_k}(e(\mu^k(i;p)))\right)\\
&= e(\mu(i;p))
\end{align*}
by the induction hypothesis since $\mu^k\in Q^n_m$. So \eqref{eq:semilatticeproducthypothesis} holds. That \eqref{eq:semilatticeproducthypothesis2} also holds is easily checked. The lemma can be proven by repeated use of this result as well as some similar computations.
\eproof

We now give a non-iterative description of the inclusion $E_k\hookrightarrow\Zz[E_{k-1}\reg]$.

\blemma\label{non-iterative-e-of-mu}
Let $\mu\in Q^n_k$. We have for $k>1$,
\begin{align*}
e(\mu)&=\\
  &\sum_{\omega\subseteq\mu_k}(-1)^{|\omega|}e(\mu^k)\prod_{p\in\omega}\bigvee\cR_p(e)\\
+ &\sum_{\substack{\rho\subseteq\mu_k \\ \rho\neq\varnothing}}\sum_{\omega\subseteq\mu_k\setminus\rho}\sum_{i=1}^{k-1}(-1)^{|\rho|+|\omega|}e(\mu^k(i;\rho))\prod_{p\in\omega}\bigvee \cR_p(e)
\end{align*}
where we take a product over the empty set to be 1.
\elemma
\bproof
We first note that if $\mu_k=\{p\}$, $\cR_p(\mu^k)=\bigcup_{i=1}^{k-1}\{e(\mu^k(i;p))\}\cup e(\mu^k)\cR_p(e)$. This follows from induction on the definition of $\cR_p(\mu^k)$ and using Lemma~\ref{semilatticeproduct-higher-level}. The same lemma gives us that $e(\mu^k(i;p))e(\mu^k(j;p))=0$ for $i\neq j$ and that $e(\mu^k(i;p))\cR_p(e)=\{0\}$ for any $i$. So
$$
e(\mu)=e(\mu^k)-\bigvee\cR_p(e(\mu^k))=e(\mu^k)-e(\mu^k)\bigvee \cR_p(e)-\sum_{i=1}^{k-1}e(\mu^k(i;p)).
$$
Then if $\mu_k$ is not a singleton we get by definition
\begin{align*}
e(\mu)&=\prod_{p\in\mu_k}\left(e(\mu^k)-e(\mu^k)\bigvee \cR_p(e)-\sum_{i=1}^{k-1}e(\mu^k(i;p))\right)\\
&= \\
&\prod_{p\in\mu_k}\left(e(\mu^k)-e(\mu^k)\bigvee \cR_p(e)\right)\\
+ &\sum_{\substack{\rho\subseteq\mu_k \\ \rho\neq\varnothing}}(-1)^{|\rho|}\prod_{p\in\mu_k\setminus\rho}\left(e(\mu^k)-e(\mu^k)\bigvee \cR_p(e)\right)\prod_{p\in\rho}\sum_{i=1}^{k-1}e(\mu^k(i;p))
\end{align*}
The first term can easily be written on the form we want it. Now by Lemma~\ref{semilatticeproduct-higher-level}, $e(\mu^k(i;p))e(\mu^k(j;p))=0$ for $i\neq j$. So we get
$$
\prod_{p\in\rho}\sum_{i=1}^{k-1}e(\mu^k(i;p))=\sum_{i=1}^{k-1}\prod_{p\in\rho}e(\mu^k(i;p))=\sum_{i=1}^{k-1}e(\mu^k(i;\rho)).
$$
Then by repeated use of Lemma~\ref{semilatticeproduct-higher-level} it is also easy to get the last term on the form we want it.
\eproof

In the rest of the section we will assume the following three properties of $\sR$.
\begin{enumerate}
\item[(A)] There is an $n\in\Nz$ such that $|\sR(e)|=n$ for all $e\in E\reg$.
\item[(B)] For each $e\in E^\times$ there is an indexing $\sR(e)=\{\cR_i(e)\}_{i=1}^n$ satisfying $g\cdot \cR_i(e)=\cR_i(g\cdot e)$ for each $1\leq i\leq n$ and $g\in\Gamma$.
\item[(C)] This indexing is also subject to the following condition. For each $i,j\in\{1,\ldots,n\}$ with $i\neq j$ there is a number $i\# j\in\{1,\ldots,n\}$ such that for each $e\in E^\times$ and $\ve\in E(\bigvee\cR_i(e))$, $\ve\cR_j(e)=\cR_{i\# j}(\ve)$.
\end{enumerate}
Roughly speaking, these properties ensure that the finite covers $\sR(e)$ exhibit a uniform behaviour with respect to the group action and multiplication. As we will see, these conditions allow us to control the behaviour the elements of $E_k$ ($k>0$) with respect to the group action of $\Gamma$, and with respect to products. This, in turn, will allow us to describe the maps $\partial_k$ in Proposition~\ref{chaincomplex-formula-C}.

It is straightforward to show that if $i,j,k\in\{1,\ldots,n\}$ are pairwise distinct, then $i\# j\neq i\# k$. For each $1\leq i\leq n$ let $M_i:\Zz[\Gamma\setminus E\reg]\to\Zz[\Gamma\setminus E\reg]$ be given by $M_i[e]=[\bigvee\cR_i(e)]$. Given $\rho\subseteq\{1,\ldots,n\}$, define $M_\rho:\Zz[\Gamma\setminus E\reg]\to\Zz[\Gamma\setminus E\reg]$ by $M_\rho[e]=\left[\prod_{p\in\rho}\left(\bigvee\cR_p(e)\right)\right]$.  Given $\rho\subseteq\{1,\ldots,n\}$ with $i\notin\rho$, let $i\#\rho\defeq\menge{i\# p}{p\in\rho}$. Given $\omega\subseteq\{1,\ldots,n\}$ and $j\notin\omega$ define iteratively $\omega\# j=(q\#(\omega\setminus\{q\}))\#(q\# j)$. One can show that this definition is independent on which order one picks out the $q$'s by using the commutativity of the semilattice product. With $\omega\cap\rho=\varnothing$, define $\omega\#\rho\defeq\menge{\omega\# p}{p\in\rho}$.

\blemma\label{M-rho-factorization}
We have for any $p\in\rho$ that $M_\rho=M_{p\#(\rho\setminus\{p\})}M_p$. 
\elemma
\bproof
Let $e\in E\reg$ and $p\in\rho$. We have that $\bigvee\cR_p(e)$ is equal to a finite sum $\sum_{\ve\in E(\bigvee\cR_p(e))}n_\ve \ve$, so
\begin{align*}
M_\rho[e]&=\left[\prod_{q\in\rho}\left(\bigvee\cR_q(e)\right)\right]=\left[\left(\bigvee\cR_p(e)\right)\prod_{q\in\rho\setminus\{p\}}\left(\bigvee\cR_q(e)\right)\right]=\\
&=\sum n_\ve \left[\ve\prod_{q\in\rho\setminus\{p\}}\left(\bigvee\cR_p(e)\right)\right]=\sum n_\ve \left[\prod_{q\in\rho\setminus\{p\}}\ve\left(\bigvee\cR_q(e)\right)\right]\\
&=\sum n_\ve \left[\prod_{q\in\rho\setminus\{p\}}\left(\bigvee\cR_{p\# q}(\ve)\right)\right]=\sum n_\ve M_{p\#(\rho\setminus\{p\})}[\ve]\\
&=M_{p\#(\rho\setminus\{p\})}\left(\sum n_\ve[\ve]\right)=M_{p\#(\rho\setminus\{p\})}M_p[e].
\end{align*}
\eproof

Given $\mu\in Q^n_k$ and $e\in E\reg$, define $[e](\mu)\defeq [e(\mu)]\in \Zz[\Gamma\setminus E_k\reg]$. This is well-defined because of property (B). Given a finite sum $x=\sum_{[e]\in\Gamma\setminus E}a_{[e]}[e]\in\Zz[\Gamma\setminus E\reg]$, define $x(\mu)\defeq \sum_{[e]\in \Gamma\setminus E}a_e[e](\mu)$. For $\mu\in Q^n_k$, $\omega\cap\mu_j=\varnothing$ for each $j$, let $\omega\# \mu\defeq \omega\#\mu_1|\cdots|\omega\#\mu_k$. Also, let $\varnothing\#\mu\defeq\mu$.

\bprop\label{chaincomplex-formula-C}
We have for $k>1$, $\mu\in Q^n_k$ and $x\in\Zz[\Gamma\setminus E\reg]$
\begin{align*}
\partial_k(x(\mu))&=\sum_{\omega\subseteq\mu_k}(-1)^{|\omega|}(M_\omega x)(\omega\# \mu^k)\\
&+\sum_{\substack{\rho\subseteq\mu_k \\ \rho\neq\varnothing}}\sum_{\omega\subseteq\mu_k\setminus\rho}\sum_{i=1}^{k-1}(-1)^{|\rho|+|\omega|}(M_\omega x)(\omega\#\mu^k(i;\rho))
\end{align*}
and for $\mu\in Q^n_1$,
$$
\partial_1(x(\mu))=\sum_{\omega\subseteq\mu_1}(-1)^{|\omega|}M_\omega x.
$$
\eprop
\bproof
It is sufficient to show that for any $k$, $\mu\in Q^n_k$, $\omega\subseteq\{1,\ldots,n\}$ with $\omega\cap\mu_j=\varnothing$ for all j, and $e\in E\reg$
\begin{equation}\label{eq:matrixformulainduction}
\left[e(\mu)\prod_{p\in\omega}\bigvee\cR_p(e)\right]=(M_\omega [e])(\omega\#\mu)
\end{equation}
If $\omega=\varnothing$, there is nothing to prove. Let $q\in\omega$. We show that $e(\mu)\ve=\ve(q\#\mu)$ for any $\ve\in E(\bigvee\cR_q(e))$. This clearly holds for $\mu\in Q^n_1$. It is then a simple induction proof to show that it holds for $\mu\in Q^n_k$ for all $k$. Assume now that we have shown that \eqref{eq:matrixformulainduction} holds for any $\omega$ of size $j$ and for any $e\in E\reg$ and $\mu\in Q^n_k$. Then if $|\omega|=j+1$, $q\in\omega$ and $\bigvee\cR_q(e)=\sum_{\ve\in E(\bigvee\cR_q(e))}n_\ve \ve$, we get by using Lemma~\ref{M-rho-factorization} that
\begin{align*}
\left[e(\mu)\prod_{p\in\omega}\bigvee\cR_q(e)\right]&=\sum n_\ve\left[e(\mu)\ve\prod_{p\in\omega\setminus\{q\}}\bigvee\cR_p(e)\ve\right]\\
&=\sum n_\ve \left[\ve(q\#\mu)\prod_{p\in\omega\setminus\{q\}}\bigvee\cR_{q\# p}(\ve)\right]\\
&=\sum n_\ve \left[\ve(q\#\mu)\prod_{p\in q\#(\omega\setminus\{q\})}\bigvee\cR_{p}(\ve)\right]\\
&=\sum n_\ve (M_{q\#(\omega\setminus\{q\})}[\ve])((q\#(\omega\setminus\{q\}))\#(q\#\mu))\\
&=\sum n_\ve (M_{q\#(\omega\setminus\{q\})}[\ve])(\omega\#\mu)\\
&=(M_{q\#(\omega\setminus\{q\})}\sum n_\ve [\ve])(\omega\#\mu)=(M_{q\#(\omega\setminus\{q\})}M_q[e])(\omega\#\mu)\\
&=(M_\omega [e])(\omega\#\mu).
\end{align*}
The result then follows by induction on $|\omega|$.
\eproof

Let $S_k$ be the permutation group on $\{1,\ldots,k\}$. Given $\sigma\in S_k$ and $\mu\in Q^n_k$, let $\sigma(\mu)=(\mu_{\sigma(1)}|\cdots|\mu_{\sigma(k)})$. Recall that an inversion for $\sigma$ is a pair $(a,b)$ with $a<b$ and $\sigma(a)>\sigma(b)$. Let $m(\sigma)$ be the number of inversions of $\sigma$ so that $(-1)^{m(\sigma)}=sgn(\sigma)=(-1)^p$, where $p$ is the number of transpositions in any given decomposition of $\sigma$ into transpositions. Let now $N^n_k\subseteq Q^n_k$ be defined by
$$
N^n_k\defeq\menge{(\mu_1|\cdots|\mu_k)}{\mu_i\in\gekl{1,\ldots,n},\mu_i<\mu_j\mbox{ when }i<j},
$$
and for $\mu\in N^n_k$ and $1\leq i\leq k$ let $\rho_{\mu,i}\in S_{k-1}$ be the unique permutation such that $\rho_{\mu,i}(\mu_i\#\mu^i)\in N^n_{k-1}$.

Define $\wt{C}$ to be the sequence
{\small
$$
0 \to \bigoplus_{\mu\in N^n_n}\Zz[\Gamma\setminus E\reg] \to \bigoplus_{\mu\in N^n_{(n-1)}}\Zz[\Gamma\setminus E\reg]\to\cdots
	\to \bigoplus_{\mu\in N^n_1}\Zz[\Gamma\setminus E\reg] \to \bigoplus_{\mu\in N^n_0}\Zz[\Gamma\setminus E\reg] \to 0
$$
}
with connecting maps $d_k:\bigoplus_{\mu\in N^n_k}\Zz[\Gamma\setminus E\reg] \to \bigoplus_{\mu\in N^n_{(k-1)}}\Zz[\Gamma\setminus E\reg]$ given by
{\small
$$
d_k(\oplus_{\mu\in N^n_k}x_\mu)=\oplus_{\lambda\in N^n_{k-1}}\sum_{\mu\in N^n_k}\sum_{i=1}^k(-1)^{i+1}\left(\delta_{\lambda,\mu^i}+(-1)^{m(\rho_{\mu,i})+1}\delta_{\lambda,\rho_{\mu,i}(\mu_i\# \mu^i)}M_{\mu_i}\right)x_\mu
$$
}
where $\delta$ is the Dirac delta. Our next goal is to show that $\wt{C}$ is a chain complex and to compare $H_*(\wt{C})$ with $H_*(C)$. This is helpful because the complexity of the first chain complex is generally much lower than that of the latter since $|Q^n_k|$ grows much faster than $k!\binom{n}{k}=k!|N^n_k|$.

For the comparison of $H_*(\wt{C})$ and $H_*(C)$, we construct a chain map $f:\wt{C}\to C$. Let $a_k$ be the integer sequence defined by $a_0=0$ and $a_k=a_{k-1}+k-1$. Define maps
$$
f_k:\bigoplus_{\mu\in N^n_k}\Zz[\Gamma\setminus E\reg]\to\Zz[\Gamma\setminus E_k\reg]
$$
by
$$
f_k(\oplus_{\mu\in N^n_k}x_\mu)=\sum_{\mu\in N^n_k}\sum_{\sigma\in S_k}(-1)^{m(\sigma)+a_k}x_\mu(\sigma(\mu))
$$
when $k\geq 1$ and let $f_0:\Zz[\Gamma\setminus E\reg]\to \Zz[\Gamma\setminus E\reg]$ be the identity.

\bprop
We have that $\wt{C}$ is a chain complex and that the homomorphisms $\{f_k\}_{k=0}^n$ give a chain map $\wt{C}\to C$.
\eprop
\bproof
We need to show that for each $1\leq k\leq n$, $\partial_k f_k=f_{k-1}d_k$. Since each $f_k$ is injective and $C$ is a chain complex it will then follow that $\tilde{C}$ is a chain complex. It is straightforward to check that $\partial_1 f_1=f_0 d_1$. Let $k>1$. Using Proposition~\ref{chaincomplex-formula-C} we have
\begin{align*}
&\partial_k f_k(\oplus_{\mu\in N^n_k}x_\mu) \\
&= \sum_{\mu\in N^n_k}\sum_{\sigma\in S_k}(-1)^{m(\sigma)+a_k}\partial_k(x_\mu(\sigma(\mu)))\\
&= \\
&\sum_{\mu\in N_k^n}\sum_{\sigma\in S_k}(-1)^{m(\sigma)+a_k}(x_\mu(\sigma(\mu)^k)-(M_{\sigma(\mu)_k}x_\mu)(\sigma(\mu)_k\#\sigma(\mu)^k))\\
- &\sum_{\mu\in N_k^n}\sum_{\sigma\in S_k}\sum_{j=1}^{k-1}(-1)^{m(\sigma)+a_k}x_\mu(\mu_{\sigma(1)}|\cdots|\mu_{\sigma(j-1)}|\mu_{\sigma(j)},\mu_{\sigma(k)}|\mu_{\sigma(j+1)}|\cdots|\mu_{\sigma(k-1)})
\end{align*}
In the last line of the previous equation, switch the two rightmost summation signs. Fix $\mu$ and $j$ and look at the expression
\begin{equation}\label{eq:permutationsum}
(-1)^{a_k}\sum_{\sigma\in S_k}(-1)^{m(\sigma)}x_\mu(\mu_{\sigma(1)}|\cdots|\mu_{\sigma(j-1)}|\mu_{\sigma(j)},\mu_{\sigma(k)}|\mu_{\sigma(j+1)}|\cdots|\mu_{\sigma(k-1)})
\end{equation}
Let $\sigma\in S_k$ with $\sigma(j)=a$ and $\sigma(k)=b$. Let $(a,b)$ be the transposition that switches $a$ and $b$. Then $\xi=(a,b)\sigma$ has opposite parity of $\sigma$, and $\xi(j)=b$, $\xi(k)=a$, and $\xi$ and $\sigma$ agree on all other elements. Now for any $x$,
\begin{align*}
&x(\mu_{\sigma(1)}|\cdots|\mu_{\sigma(j-1)}|\mu_{\sigma(j)},\mu_{\sigma(k)}|\mu_{\sigma(j+1)}|\cdots|\mu_{\sigma(k-1)})=\\
&x(\mu_{\xi(1)}|\cdots|\mu_{\xi(j-1)}|\mu_{\xi(j)},\mu_{\xi(k)}|\mu_{\xi(j+1)}|\cdots|\mu_{\xi(k-1)}).
\end{align*}
Moreover, $(a,b)\xi=\sigma$. Since every $\sigma\in S_k$ comes paired with such a $\xi$ of opposite parity (i.e. $(-1)^{m(\sigma)}=-(-1)^{m(\xi)}$) we get that expression \eqref{eq:permutationsum} is equal to $0$ for all $\mu$ and $j$. We are left with
$$
\partial_k f_k(\oplus_{\mu\in N^n_k}x_\mu)=\sum_{\mu\in N_k^n}\sum_{\sigma\in S_k}(-1)^{m(\sigma)+a_k}(x_\mu(\sigma(\mu)^k)-(M_{\sigma(\mu)_k}x_\mu)(\sigma(\mu)_k\# \sigma(\mu)^k))
$$
On the other hand, we have
\begin{align*}
&f_{k-1}d_k(\oplus_{\mu\in N^n_k}x_\mu)\\
&=\sum_{\lambda\in N^n_{k-1}}\sum_{\sigma\in S_{k-1}}\sum_{\mu\in N^n_k}\sum_{i=1}^k (-1)^{i+1}(-1)^{m(\sigma)+a_{k-1}}\\
&\cdot(\delta_{\lambda,\mu^i}x_\mu(\sigma(\lambda))+(-1)^{m(\rho_{\mu,i})+1}\delta_{\lambda,\rho_{\mu,i}(\mu_i\# \mu^i)}(M_{\mu_i}x_\mu)(\sigma(\lambda)))\\
&=\sum_{\mu\in N^n_k}\sum_{i=1}^k\sum_{\sigma\in S_{k-1}}(-1)^{m(\sigma)+a_{k-1}+i+1}\\
&\cdot(x_\mu(\sigma(\mu^i))+(-1)^{m(\rho_{\mu,i})+1}(M_{\mu_i}x_\mu)((\sigma\rho_{\mu,i})(\mu_i\# \mu^i)))\\
&=\sum_{\mu\in N^n_k}\sum_{i=1}^k\sum_{\sigma\in S_{k-1}}(-1)^{m(\sigma)+a_{k-1}+i+1}(x_\mu(\sigma(\mu^i))-(M_{\mu_i}x_\mu)(\sigma(\mu_i\# \mu^i))).
\end{align*}
The last equality holds since it does not matter whether we sum over $S_{k-1}$ or $S_{k-1}\rho_{\mu,i}$. We also have for each $1\leq j\leq n$,
$$
(\sigma(\mu_i\#\mu^i))_j=(\mu_i\#\mu^i)_{\sigma(j)}=\mu_i\#(\mu^i)_{\sigma(j)}=\mu_i\#(\sigma(\mu^i))_j=(\mu_i\#\sigma(\mu^i))_j
$$
so $\sigma(\mu_i\#\mu^i)=\mu_i\#\sigma(\mu^i)$. It remains to show that for each $\mu\in N^n_k$,
\begin{align*}
&\sum_{\xi\in S_k}(-1)^{m(\xi)+a_k}(x_\mu(\xi(\mu)^k)-(M_{\mu_{\xi(k)}}x_\mu)(\mu_{\xi(k)}\#(\xi(\mu)^k)))\\
=&\sum_{\sigma\in S_{k-1}}\sum_{i=1}^k(-1)^{m(\sigma)+a_{k-1}+i+1}(x_\mu(\sigma(\mu^i))-(M_{\mu_i}x_\mu)(\mu_i\# \sigma(\mu^i)))
\end{align*}
It is sufficient to find a bijection $\xi:S_{k-1}\times\{1,\ldots,k\}\to S_k$ satisfying that for all $\sigma,i$, $\xi(\sigma,i)(k)=i$, $(\xi(\sigma,i)(\mu))^k=\sigma(\mu^i)$ and $m(\xi(\sigma,i))+a_k\equiv m(\sigma)+a_{k-1}+i+1 \mbox{(mod 2)}$. Define
\[
\xi(\sigma,i)(j)=\begin{cases}
i&\mbox{ if }j=k\\
\sigma(j)&\mbox{ if }\sigma(j)<i\\
\sigma(j)+1&\mbox{ if }i\leq\sigma(j)
\end{cases}
\]
It is straightforward to check that $\xi$ satisfies the right properties.
\eproof

\bremark
We have been able to verify that $f_*$ induces an isomorphism $H_*(\wt{C})\cong H_*(C)$ when $n=2$. We have also verified that it induces an isomorphism for $n=3$ when $i\#j=j$ for all $i\neq j$. We suspect that $f_*$ always induces an isomorphism, but we have not been able to find a general proof.
\eremark

\section{Examples}

We now present an application of the ideas we developed in the previous sections to group homology and cohomology. More precisely, for certain groups $\Gamma$, we obtain free $\Zz \Gamma$-resolutions of the trivial $\Zz \Gamma$-module $\Zz$. Such resolutions play a fundamental role in group homology and cohomology. Their existence is guaranteed, but only in an abstract sense. Typically, concrete constructions of such resolutions allow for precise computations and a detailed analysis of the homological behaviour of the groups of interest. For instance, in our case, we obtain resolutions of finite length. This leads to upper bounds on the cohomological dimension. We refer to \cite{Bro} for more information about group homology and cohomology. In the following, we first follow a general approach using the language of subsemigroups of groups. Afterwards, we give concrete examples of groups given by generators and relations where our ideas work.

\subsection{Algebraic independent resolutions for $\Zz$, and group (co)homology}
\label{groups}

Let $P$ be a unital subsemigroup of a group $G$, let $\Sigma$ be a set of generators for $P$, and let us assume that the following conditions are fulfilled:
\begin{enumerate}
\item[(a)] $P \subseteq G$ is left quasi-lattice ordered, $P$ is right Ore and its enveloping group is $G$.
\item[(b)] For all $x$, $y$ and $z$ in $G$ with $xP \cap yP = zP$ and for every $s \in \Sigma$, either $xsP \cap yP = zP$ or there exists $t \in \Sigma$ with $xsP \cap yP = ztP$.
\item[(c)] For every finite subset $F$ of $\Sigma$ and $s \in \Sigma \setminus F$, we have $sP \cap \bigcap_{t \in F} tP \subsetneq \bigcap_{t \in F} tP$.
\end{enumerate}
(a) means that $P$ contains no non-trivial units, and that for every $x$ and $y$ in $G$, either $xP \cap yP = \varnothing$ or there exists $z \in G$ with $xP \cap yP = zP$. (a) implies that $E \defeq \gekl{gP}_{g \in G} \cup \gekl{\varnothing}$ is a semilattice with respect to intersection. Consider the collection $\sR$ of finite covers for $E$ defined by $\sR(gP) \defeq \menge{\gekl{gsP}}{s \in \Sigma}$. $\gekl{gsP}$ is a cover of $gP$ because $P$ is right Ore. It turns out that $G \curvearrowright E$ and $\sR$ satisfy conditions (i) to (iii). (iii) is obvious, (ii) follows from condition (c) and (i) follows from (b). Moreover, let $I = \spkl{\menge{e - e(\cR)}{e \in E\reg, \ \cR \in \sR(e)}}_{\Zz} \triangleleft \Zz[E\reg]$. In the quotient $\Zz[E\reg] / I$, we identify $gP$ with $gsP$ for every $g \in G$ and $s \in \Sigma$. Therefore, all the idempotents in $E\reg$ are identified, and we end up with $\Zz[E\reg] / I \cong \Zz$. 

Therefore, by Theorem~\ref{maintheo}, we obtain semilattices $E$, $E_1$, $E_2$, ..., together with $G$-actions, and a $G$-equivariant long exact sequence
\bgl
\label{free-resolution--Z}
  \dotso \to \Zz[E_2\reg] \to \Zz[E_1\reg] \to \Zz[E\reg] \to \Zz \to 0.
\egl
$G$-equivariance just means that this is a sequence of $\Zz G$-modules. Even more, since $G$ acts freely on $E\reg$, Corollary~\ref{cor-stabilizer} implies that $G$ acts freely on $E_n\reg$ for every $n \in \Nz$. In other words, \eqref{free-resolution--Z} is a free resolution of $\Zz$ over $\Zz G$. And if $\abs{\Sigma} < \infty$, then \eqref{free-resolution--Z} becomes a free resolution of $\Zz$ over $\Zz G$ of length at most $\abs{\Sigma}$. We obtain the following
\bcor
If a group $G$ has a subsemigroup $P$ satisfying (a), (b) and (c) with $\abs{\Sigma} < \infty$, then $G$ is of type FL, and the cohomological dimension of $G$ is at most $\abs{\Sigma}$. In particular, $H_i(G,M) = \gekl{0}$ and $H^i(G,M) = \gekl{0}$ for every $i > \abs{\Sigma}$ and for every $\Zz G$-module $M$. 
\ecor
The reader may consult \cite{Bro} for more details about group (co)homology. In particular, cohomological dimension is explained in \cite[VIII, \S~2]{Bro}, and type FL is explained in \cite[VIII, \S~6]{Bro}.

\subsection{Concrete examples of groups}

Let $\Sigma$ be a finite set, $\sigma$: $\Sigma \times \Sigma \to \Sigma \times \Sigma$ a bijection. We write the first component of $\sigma$ as $\sigma_l$ and the second one as $\sigma_r$, so that $\sigma(a,b) = (\sigma_l(a,b),\sigma_r(a,b))$. Assume further that $\sigma$ satisfies the following conditions:
\begin{itemize}
\item[($*$)] $\sigma \vert_{\Delta} = \id_{\Delta}$, where $\Delta \subseteq \Sigma$ is the diagonal $\menge{(a,a)}{a \in \Sigma}$.
\item[($*$$*$)] $\sigma$ is flip-invariant: $\sigma(b,a) = (\sigma_r(a,b),\sigma_l(a,b))$.
\item[($*$$*$$*$)] For fixed $a \in \Sigma$, $\Sigma \setminus \gekl{a} \ni x \ma \sigma_l(a,x) \in \Sigma$ is injective.
\item[($*$$*$$*$$*$)] For pairwise distinct $a$, $b$ and $c$ in $\Sigma$, let $\sigma(a,b) = (d,e)$, $\sigma(b,c) = (f,g)$ and $\sigma(c,a) = (h,i)$. Then there exist $j$, $k$ and $l$ in $\Sigma$ with $\sigma(e,f) = (j,k)$, $\sigma(g,h) = (k,l)$ and $\sigma(i,d) = (l,j)$. In addition, $\sigma^{-1}$ has this property as well.
\end{itemize}

\blemma
Assume that $\sigma$ is as above satisfying ($*$) to ($*$$*$$*$$*$). If we now let $P$ be the (unital) semigroup $\spkl{\Sigma,R}^+$ generated by $\Sigma$ subject to the relations $R$ given by $a \sigma_l(a,b) = b \sigma_r(a,b)$, $a$ and $b$ in $\Sigma$, and if $G$ is the corresponding group $\spkl{\Sigma,R}$, then $P$ is (in a canonical way) a subsemigroup of $G$, and $P$, $G$ and $\Sigma$ satisfy (a), (b) and (c) from \S~\ref{groups}.
\elemma
\bproof
We first prove that $P$ embeds into $G$ and that (a) holds for $P \subseteq G$. First of all, $(\Sigma,R)$ is $r$-complete and $l$-complete in the sense of \cite{Deh}. We show that ($*$), ($*$$*$) and ($*$$*$$*$$*$) imply that $(\Sigma,R)$ is $r$-complete; $l$-completeness then follows by symmetry. Since the relations in $R$ preserve the length of words in $\Sigma$, it is clear that $(\Sigma,R)$ is $r$-homogeneous. We now show that ($*$$*$$*$$*$) implies that the strong $r$-cube condition is satisfied on $\Sigma$: Let $a$, $b$ and $c$ be pairwise distinct elements in $\Sigma$ such that $\sigma(a,b) = (d,e)$, $\sigma(b,c) = (f,g)$ and $\sigma(c,a) = (h,i)$. By ($*$$*$$*$$*$), there exist $j$, $k$ and $l$ in $\Sigma$ with $\sigma(e,f) = (j,k)$, $\sigma(g,h) = (k,l)$ and $\sigma(i,d) = (l,j)$. Therefore, using the notation from \cite{Deh}, we have $a^{-1} b b^{-1} c \curvearrowright_r d e^{-1} f g^{-1} \curvearrowright_r d j k^{-1} g^{-1}$. To verify the strong $r$-cube condition, we have to show $(dj)^{-1} a^{-1} c (gk) \curvearrowright_r \ve$. Indeed, $(dj)^{-1} a^{-1} c (gk) = j^{-1} d^{-1} a^{-1} c g k \curvearrowright_r j^{-1} d^{-1} i h^{-1} g k \curvearrowright_r j^{-1} j l^{-1} l k^{-1} k \curvearrowright_r \ve$. Hence $(\Sigma,R)$ is $r$-complete.

\cite[Corollary~6.2]{Deh} implies that $P = \spkl{\Sigma,R}^+$ is cancellative because $\sigma$ is bijective with $\sigma \vert_{\Delta} = \id_{\Delta}$.

As $\sigma$ is bijective, \cite[Proposition~6.7]{Deh} and \cite[Remark~6.9]{Deh} imply that common right multiples exist in $P$, i.e. for all $p$ and $q$ in $P$, $pP \cap qP \neq \varnothing$. By symmetry, common left multiplies in $P$ exist as well. In other words, $P$ is both left and right Ore. Of course, the enveloping group (or group of fractions) of $P$ is $G$.

In addition, \cite[Proposition~6.10]{Deh} yields that $P$ admits least common right multiples, i.e., for all $p$ and $q$ in $P$, there exists $r$ with $pP \cap qP = rP$. This again follows from our conditions that $\sigma$ is bijective and $\sigma \vert_{\Delta} = \id_{\Delta}$. By symmetry, $P$ also admits least common left multiples.

Obviously, $P$ contains no non-trivial units.

All in all, it follows that $P$ embeds into $G$ and that $P$ and $G$ satisfy (a) from \S~\ref{groups}. 

For (b), we first argue that it suffices to consider the case where $x$, $y$ and $z$ are in $P$. Namely, given $x$, $y$ and $z$ in $G$ as in (b), we can find $p \in P$ with $px, py, pz \in P$ as $P$ is left Ore. So $pxP \cap pyP = pzP$. If (b) holds for elements in $P$, then for arbitrary $s \in \Sigma$, we either have $pxsP \cap pyP = pzP$ or there exists $t \in \Sigma$ with $pxsP \cap pyP = pztP$. Now multiply with $p^{-1}$ from the left, and we are done. It remains to prove (b) for elements in $P$. Assume $x$ and $y$ are in $P$. Write $x = x_1 \dotsm x_l$ and $y = y_1 \dotsm y_m$ with $x_j$ and $y_k$ in $\Sigma$. Then we can find $a_h$ and $b_i$ in $\Sigma$ with $pP \cap qP = x_1 \dotsm x_l a_1 \dotsm a_{\lambda} P = y_1 \dotsm y_m b_1 \dotsm b_{\mu} P$. So $z$ in (b) is $x_1 \dotsm x_l a_1 \dotsm a_{\lambda} = y_1 \dotsm y_m b_1 \dotsm b_{\mu}$. \cite[Proposition~6.10]{Deh} tells us that this is equivalent to $y^{-1} x = y_l^{-1} \dotsm y_m^{-1} x_1 \dotsm x_l \curvearrowright_r b_1 \dotsm b_{\mu} a_{\lambda}^{-1} \dotsm a_1^{-1}$. Therefore, given $s \in \Sigma$, it follows that $y^{-1} xs = y_l^{-1} \dotsm y_m^{-1} x_1 \dotsm x_l s \curvearrowright_r b_1 \dotsm b_{\mu} a_{\lambda}^{-1} \dotsm a_1^{-1} s$. Now, there are two possiblities: Either there exists $t$ and $c_1$, ..., $c_{\nu}$ in $\Sigma$ with $a_{\lambda}^{-1} \dotsm a_1^{-1} s \curvearrowright_r t c_1^{-1} \dotsm c_{\nu}^{-1}$, so that $xsP \cap yP = y_1 \dotsm y_m b_1 \dotsm b_{\mu} t P = ztP$, or there are $c_1$, ..., $c_{\nu}$ in $\Sigma$ with $a_{\lambda}^{-1} \dotsm a_1^{-1} s \curvearrowright_r c_1^{-1} \dotsm c_{\nu}^{-1}$, so that $xsP \cap yP = y_1 \dotsm y_m b_1 \dotsm b_{\mu} P = zP$. This proves (b).

To prove (c), let $F$ be a finite subset of $\Sigma$, say $F = \gekl{s_1, \dotsc, s_n}$. Furthermore, take $s \in \Sigma \setminus F$. We have to show that $\bigcap_{i=1}^n s_i P \neq sP \cap \bigcap_{i=1}^n s_i P$. We proceed inductively on $n$. For $n=1$, i.e. $F = \gekl{t}$, then if $\sigma(s,t) = (a,b)$, the relation $sa = tb$ implies that $sP \cap tP = saP = tbP$ by \cite[Proposition~6.10]{Deh}. Now assume that our assertion is true for every $F' \subseteq \Sigma$ with $\abs{F'} < \abs{F}$. It is straightforward to check that there exist $x_1$, ..., $x_n$ in $\Sigma$ such that for every $1 \leq m \leq n$, we have $s_1 P \cap \dotso \cap s_m P = x_1 \dotsm x_m P$. Therefore
\bgl
\label{contra1}
  s_n^{-1} x_1 \dotsm x_{n-2} \curvearrowright_r y_1 \dotsm y_{n-2} z^{-1}
\egl
and
\bgl
\label{contra2}
  z^{-1} x_{n-1} \curvearrowright_r y_{n-1} x_n^{-1},
\egl
for some $y_1$, ..., $y_{n-1}$ and $z$ in $\Sigma$. As a consequence of \eqref{contra1}, we get
\bgl
\label{contra3}
  s_1 P \cap \dotso \cap s_{n-2} P \cap s_n P = x_1 \dotsm x_{n-2} z P.
\egl
Also, using the induction hypothesis, $s^{-1} x_1 \dotsm x_{n-2} \curvearrowright_r a_1 \dotsm a_{n-2} b^{-1}$, $b^{-1} x_{n-1} \curvearrowright_r a_{n-1} c^{-1}$, so that $s^{-1} x_1 \dotsm x_{n-1} \curvearrowright_r a_1 \dotsm a_{n-1} c^{-1}$, for some $a_1$, ..., $a_{n-1}$, $b$ and $c$ in $\Sigma$. If $c$ were equal to $x_n$, then
\bgl
\label{contra4}
  x_{n-1} x_n = b a_{n-1}
\egl
would be one of the relations in $R$. But $b$ cannot be equal to $z$. Otherwise, we would have $s^{-1} x_1 \dotsm x_{n-2} z \curvearrowright_r a_1 \dotsm a_{n-2} z^{-1} z \curvearrowright_r a_1 \dotsm a_{n-2}$, so that $sP \cap (x_1 \dotsm x_{n-2} z P) = x_1 \dotsm x_{n-2} z P$. But this, together with \eqref{contra}, contradicts the induction hypothesis. So $b \neq z$. However, \eqref{contra2} means that $x_{n-1} x_n = z y_{n-1}$ is a relation in $R$, and we also have relation \eqref{contra4}. Since $b \neq z$, this is a contradiction because by construction, there is exactly one relation of the form $x_{n-1} x = z y$ (for some $x$ and $y$ in $\Sigma$). So we conclude that $c \neq x_n$. Thus $s^{-1} x_1 \dotsm x_n \curvearrowright_r a_1 \dotsm a_{n-1} c^{-1} x_n \curvearrowright a_1 \dotsm a_{n-1} a_n d^{-1}$ for some $a_n$ and $d$ in $\Sigma$. This means that $sP \cap s_1 P \cap \dotso \cap s_n P = sP \cap (x_1 \dotsm x_n P) = x_1 \dotsm x_n d P \neq s_1 P \cap \dotso \cap s_n P$, and we are done.

All in all, we have proven (a), (b) and (c) for $P = \spkl{\Sigma,R}^+$, $G = \spkl{\Sigma,R}$ and the generators $\Sigma$.
\eproof

Hence, in the situation of the previous lemma, \S~\ref{groups} gives us a way to compute group homology for $G$. 

\bexample
If $\abs{\Sigma} = 2$, say $\Sigma = \gekl{a,b}$, then there are two possibilities for $\sigma$: Either $\sigma(a,b) = (b,a)$ or $\sigma(a,b) = (a,b)$. In the first case, $G = \spkl{ a,b \ \vert \ ab = ba } \cong \Zz^2$, and our method gives $H_0(G) \cong \Zz$, $H_1(G) \cong \Zz^2$, $H_2(G) \cong \Zz$, and $H_n(G) \cong \gekl{0}$ for $n>2$, as expected. In the second case, $G = \spkl{ a,b \ \vert \ a^2 = b^2 }$, and our method gives $H_0(G) \cong \Zz$, $H_1(G) \cong \Zz \oplus \Zz / 2 \Zz$, $H_2(G) \cong \gekl{0}$, and $H_n(G) \cong \gekl{0}$ for $n>2$.
\eexample

\bquestion
There are some obvious questions which come to mind about these groups $G = \spkl{\Sigma,R}$ we introduced. For instance, are they amenable? Do they satisfy the Baum-Connes conjecture? How fast does the number of isomorphism classes of these groups grow with $\abs{\Sigma}$? Is there a different method to compute group homology, maybe of a more topological nature, using nice models for classifying spaces?
\equestion

\end{document}